%% file: PCF.tex
\pdfoutput=1

\documentclass[11pt,a4paper,twoside]{article}
\usepackage{amsmath,amsxtra,amsthm,amssymb,amsfonts,wasysym,
bbm,pifont,color,subfig,wrapfig,epsfig,multirow,enumerate,tensor,
epstopdf,graphicx,cite,mathrsfs}
\usepackage[all]{xy}
\usepackage[margin=2cm,inner=2cm]{geometry}
\usepackage[bottom]{footmisc}

\newtheorem*{theorem'}{Theorem}
\newtheorem*{lemma'}{Lemma}
\newtheorem*{proposition'}{Proposition}
\newtheorem*{conjecture'}{Conjecture}
\newtheorem*{claim'}{Claim}
\newtheorem*{overview'}{Overview}

\newtheorem{theorem}{Theorem}[section]
\newtheorem{lemma}[theorem]{Lemma}
\newtheorem{proposition}[theorem]{Proposition}
\newtheorem{corollary}[theorem]{Corollary}

\newtheorem{conjecture}[theorem]{Conjecture}
\newtheorem{hypothesis}[theorem]{Hypothesis}
\theoremstyle{definition}
\newtheorem{definition}[theorem]{Definition}
\newtheorem{algorithm}[theorem]{Algorithm}
\newtheorem{remark}[theorem]{Remark}

\newcommand{\dd}{\mathrm{d}}
\newcommand{\ee}{\mathrm{e}}
\newcommand{\ZZ}{\mathbbm{Z}}
\newcommand{\LL}{\mathbbm{L}}
\newcommand{\GG}{\mathbf{G}}
\newcommand{\La}{\mathbf{\Lambda}}
\newcommand{\pr}{\mathbbm{P}}
\newcommand{\ex}{\mathbbm{E}}

\newcommand{\bp}{\vspace{-0.5\baselineskip}\begin{proof}}
\newcommand{\bpc}{\vspace{-0.5\baselineskip}\begin{proof}[Proof of Claim]}
\newcommand{\ep}{\renewcommand{\qedsymbol}{$\square$}\end{proof}}
\newcommand{\epc}{\renewcommand{\qedsymbol}{$\square$}\end{proof}}

\numberwithin{equation}{section}
\numberwithin{figure}{section}
\numberwithin{table}{section}
\begin{document}

\title{Percolation with constant freezing}
\author{Edward Mottram
\footnote{University of Cambridge, Centre for Mathematical Sciences, Wilberforce Road, Cambridge. CB3 0WA. UK. $_{}$ $_{}$ $_{}$ $_{}$ $_{}$ $_{}$ $_{}$ $_{}$
E-mail: e.j.mottram@maths.cam.ac.uk}}
\date{\today}
\maketitle

\begin{abstract}
\vspace{-2\baselineskip} \[\]
We introduce and study a model of \emph{percolation with constant freezing} (\emph{PCF}) where edges open at constant rate $1$, and clusters freeze at rate $\alpha$ independently of their size.  Our main result is that the infinite volume process can be constructed on any amenable vertex transitive graph.  This is in sharp contrast to models of percolation with freezing previously introduced, where the limit is known not to exist.  Our interest is in the study of the percolative properties of the final configuration as a function of $\alpha$.  We also obtain more precise results in the case of trees.  Surprisingly the algebraic  exponent for the cluster size depends on the degree, suggesting that there is no lower critical dimension for the model.  Moreover, even for $\alpha<\alpha_c$, it is shown that finite clusters have algebraic tail decay, which is a signature of self organised criticality.
Partial results are obtained on $\ZZ^d$, and many open questions are discussed.
\end{abstract}

\setlength{\parindent}{0pt}
\setlength{\parskip}{0.5\baselineskip}


\paragraph{AMS Subject Classification:}
Primary: 60K35; Secondary: 82C43.

\input{PCF_intro}

\input{PCF_construction}

\input{PCF_Zd}

\input{PCF_tree}

\input{PCF_simulations}

\setlength{\parskip}{0.5\baselineskip}

\appendix
\section*{Acknowledgements}

I would like to thank Nathana\"el Berestycki for suggesting this problem to me, for helpful discussion and for careful reading of various drafts of this paper.  I am also grateful to Geoffrey Grimmett for suggesting I use ergodic decomposition in a similar way to \cite{gkn} -- see Lemma \ref{lem:decomposition}.

This work has been supported by the UK Engineering and Physical Sciences Research Council (EPSRC) grant EP/H023348/1.

\bibliography{PCF_refs}{}
\bibliographystyle{plain}

\end{document}

%% file: PCF_intro.tex
\section{Introduction}\label{sec:intro}

Let $\alpha>0$, and let $\GG=(\mathcal{V},\mathcal{E})$ be a finite graph for now.  We consider a modification of the percolation process defined as follows.  Each edge $e\in\mathcal{E}$ opens independently at rate $1$, and each open cluster freezes independently at a constant rate $\alpha$ (regardless of its size).  Once a cluster has frozen all its neighbouring edges will remain closed forever.  We are interested in the final configuration of the edges of $\GG$, and its dependence on $\alpha$.  This is the \emph{percolation with constant freezing} (\emph{PCF}) model, and was introduced by Ben--Naim and Krapivski in the mean field case -- where several interesting features were shown (through not entirely rigorous methods) \cite{bn&k2}.  In this paper we prove that an infinite volume process can be defined on the finite dimensional lattice $\ZZ^d$ -- or more generally on any amenable vertex transitive graph $\GG$.  The existence of this process on any countable tree $\mathbf{T}$ is also shown.

\subsection{Existence of an infinite volume limit}

It is straightforward to construct the PCF model on any finite graph $\GG=(\mathcal{V},\mathcal{E})$ -- see Definition \ref{def:PCFonfinite}.  By running the process until all the clusters are frozen, we induce the \emph{PCF measure}  $\mu_{\mathbf{G},\alpha}$ on $\{0,1\}^{\mathcal{E}}$.  Our aim is to construct a PCF measure when $\GG$ is an amenable vertex transitive graph.

The infinite volume version of PCF can be understood in terms of \emph{local limits}.  Given $\GG$, we let $\GG_1\subseteq\GG_2\subseteq\ldots$ be a fixed exhaustion of $\GG$ with each $\GG_n$ finite.  Then for any finite $\mathbf{\Lambda\subseteq G}$, there exists an $N$ such that $\mathbf{\Lambda\subseteq G}_n$ for all $n\geq N$.  It is shown that the restriction of PCF on $\GG_n$ to $\mathbf{\Lambda}$ tends in law to a unique limiting process as $n\longrightarrow\infty$.  This limiting process is then equal in law to the infinite volume version of PCF restricted to $\mathbf{\Lambda}$.

\begin{theorem}\label{thm:amenable}
For every amenable vertex transitive graph $\GG=(\mathcal{V},\mathcal{E})$, and every fixed rate of freezing $\alpha > 0$, there exists an infinite volume PCF process on $\GG$ in the sense of local limits.  This induces the PCF measure $\mu_{\GG,\alpha}$ on $\{0,1\}^{\mathcal{E}}$.  Furthermore, the measure $\mu_{\GG,\alpha}$ is translation invariant.
\end{theorem}

This result is in sharp contrast to Aldous's \emph{Frozen percolation} model \cite{aldous}, in which a clusters freezes as soon as it becomes infinite.  Although Aldous showed that this model can be rigorously defined on a binary tree -- where it exhibits interesting behaviour -- it can not be constructed on a general graph.  In particular Benjamini and Schramm proved that there can be no such model model on $\ZZ^2$, \cite{b&s}.  See \cite{infinitefreeze1} for details; and see \cite{vdb&n1} and \cite{vdb&n2} for additional interesting work relating to frozen percolation. 

The construction of an infinite volume PCF process goes via a secondary process which we call \emph{warm PCF}.  Here the definition of the process is adapted so that clusters on the boundary of $\GG_n\subseteq\GG$ do not freeze.  By modifying the process in this way we obtain a form of monotonicity which in turn leads to the infinite volume limit.  We note that the technique of using a warm boundary can also be adapted to the construction of similar models.  In particular the monotonicity properties can be used to construct an infinite volume version of the Drossel Schwabl forest fire model, \cite{drossel}. 
Indeed, the construction of a stationary measure for the forest fire model in \cite{stahl} uses similar ideas.

\subsection{Properties of PCF}

\paragraph{On a lattice}

Having seen that PCF can be constructed on any amenable vertex transitive graph it is natural to ask what the process looks like.  Clearly the rate of freezing controls the behaviour of the model, and so we make the following conjecture.

\begin{conjecture}
For every dimension $d\geq 2$ there is a critical value $\alpha_c$ such that when we run rate $\alpha$ PCF on $\ZZ^d$ then if $\alpha>\alpha_c$ all clusters are almost surely finite (the \emph{sub-critical regime}), and if $\alpha<\alpha_c$ the final distribution contains infinite clusters (the \emph{super-critical regime}).
\end{conjecture}

The following proposition is a first step in this direction.

\begin{proposition}\label{prop:no_infinite}
Suppose we run PCF on the $d$ dimensional lattice $\mathbbm{Z}^d$, then provided $\alpha > 0$ is sufficiently large (depending on $d$) all clusters will remain finite almost surely.
\end{proposition}

\begin{figure}[p]
\centering
\setlength\fboxsep{0pt}
\setlength\fboxrule{1pt}
\fbox{\includegraphics[width=0.9\textwidth]{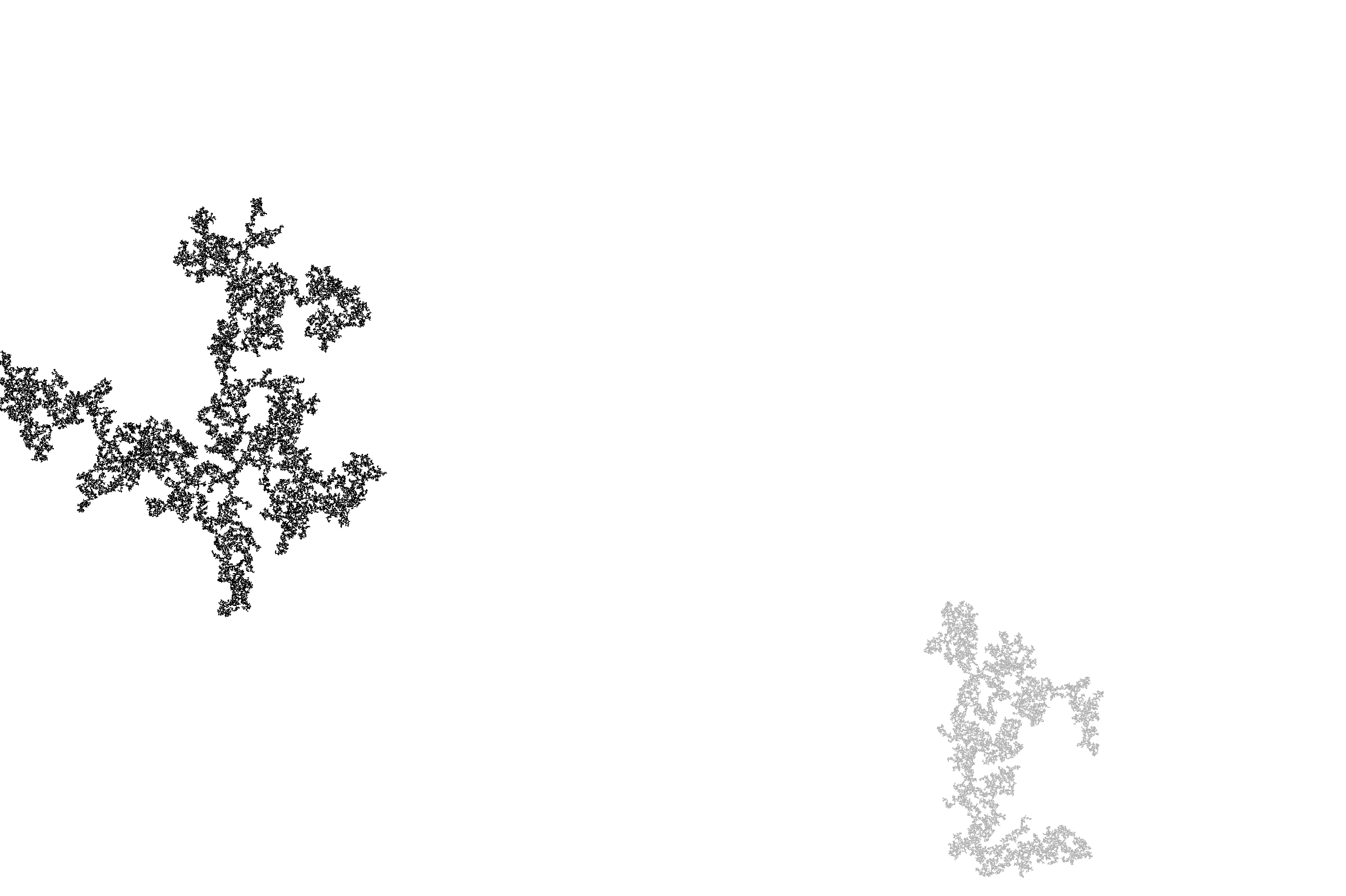}}
\caption{Sub-critical PCF \\ {The largest two clusters generated in a simulation of sub-critical PCF on a 3076 by 2048 square grid.  Here we have $\alpha=0.60$ which our simulations suggest is sufficiently quick to prevent the formation of a spanning cluster -- see Section \ref{sec:sim} for details.}}
\label{fig:PCF060}
\end{figure}

\begin{figure}[p]
\centering
\setlength\fboxsep{0pt}
\setlength\fboxrule{1pt}
\fbox{\includegraphics[width=0.9\textwidth]{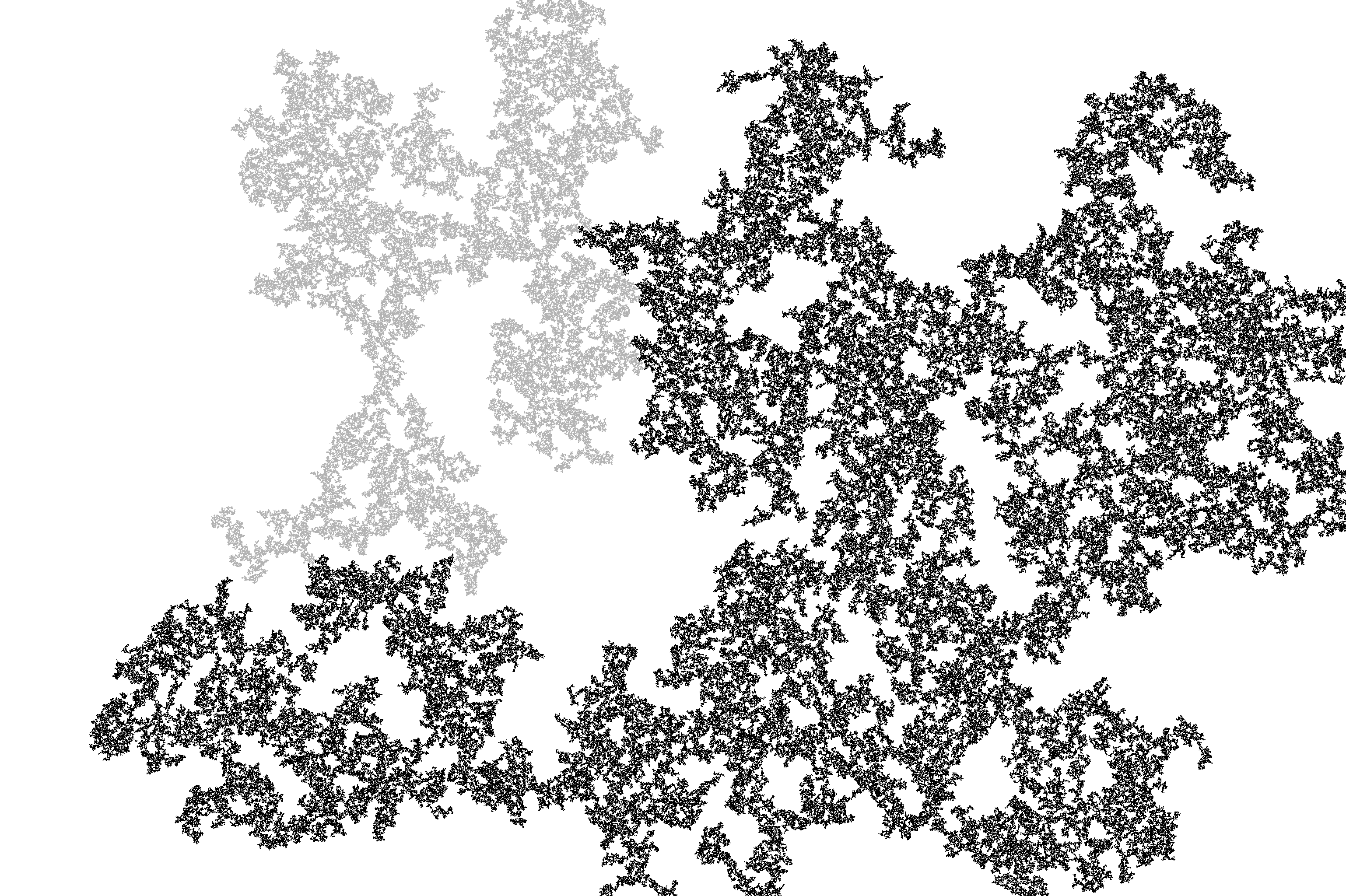}}
\caption{Near-critical PCF \\ {Our simulations suggest that on $\ZZ^2$ we have $\alpha_c \approx 0.55$.  Here are the two largest clusters generated in a simulation of PCF on a 3076 by 2048 square grid at this near critical value.}}
\label{fig:PCF055}
\end{figure}

\begin{figure}[p]
\centering
\setlength\fboxsep{0pt}
\setlength\fboxrule{1pt}
\fbox{\includegraphics[width=0.9\textwidth]{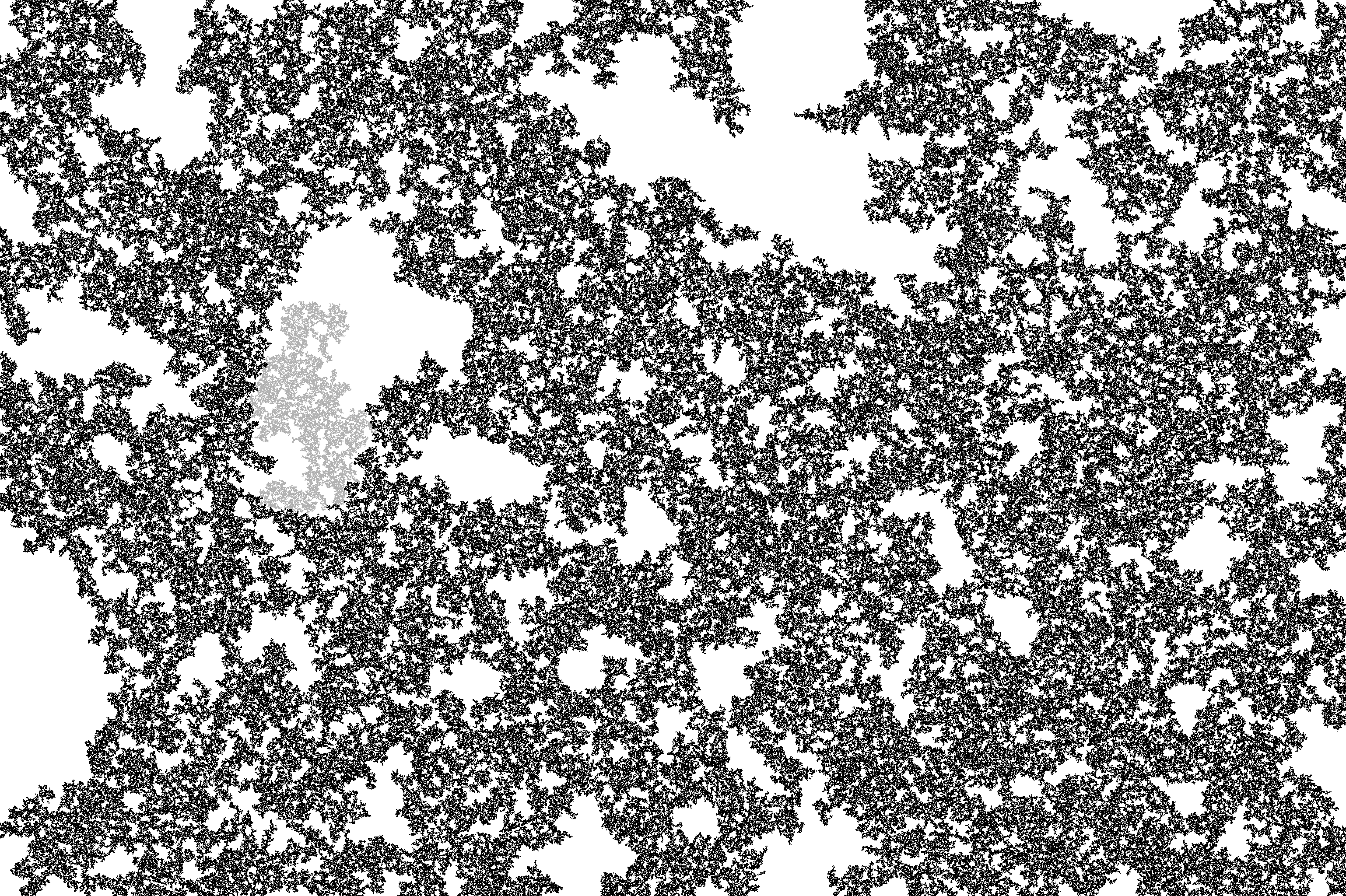}}
\caption{Super-critical PCF \\ {The largest two clusters generated in a simulation of super-critical PCF on a 3076 by 2048 square grid.  Here we have $\alpha=0.50$, which appears to be too small to prevent the formation of a spanning cluster.  In contrast to the percolation model, simulations also suggest that the size of the finite components (and thus the size of the holes) have a power law distribution.}}
\label{fig:PCF050}
\end{figure}

\begin{figure}[p]
\centering
\includegraphics[width=0.45\textwidth]{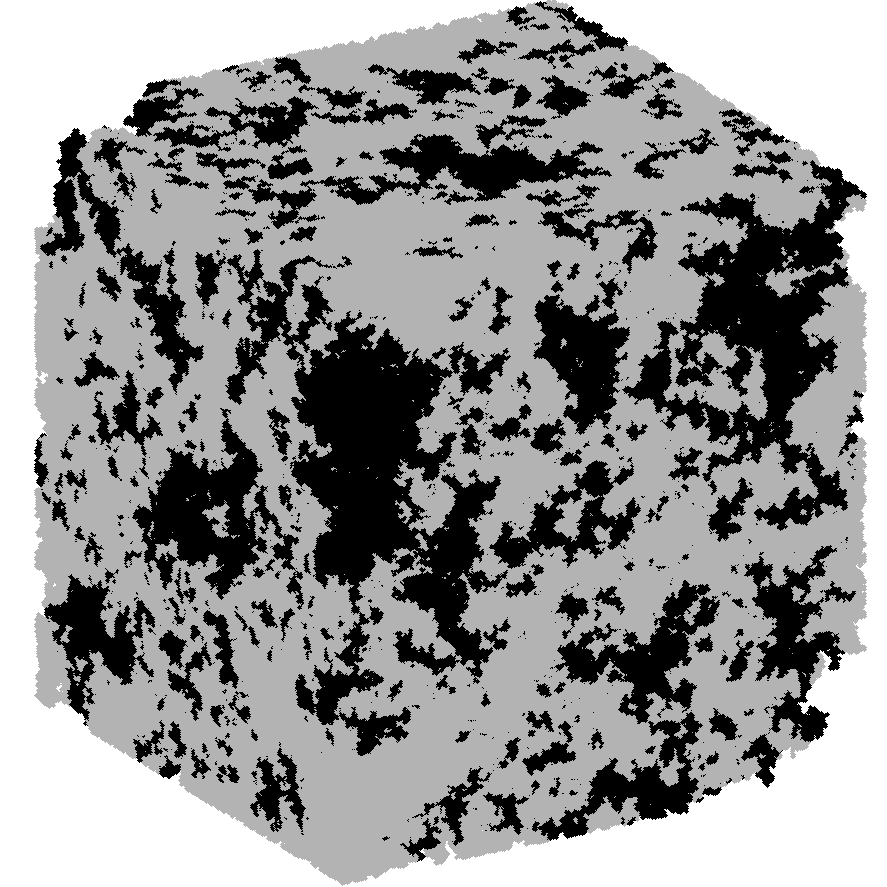}
\hspace{0.00\textwidth}
\includegraphics[width=0.45\textwidth]{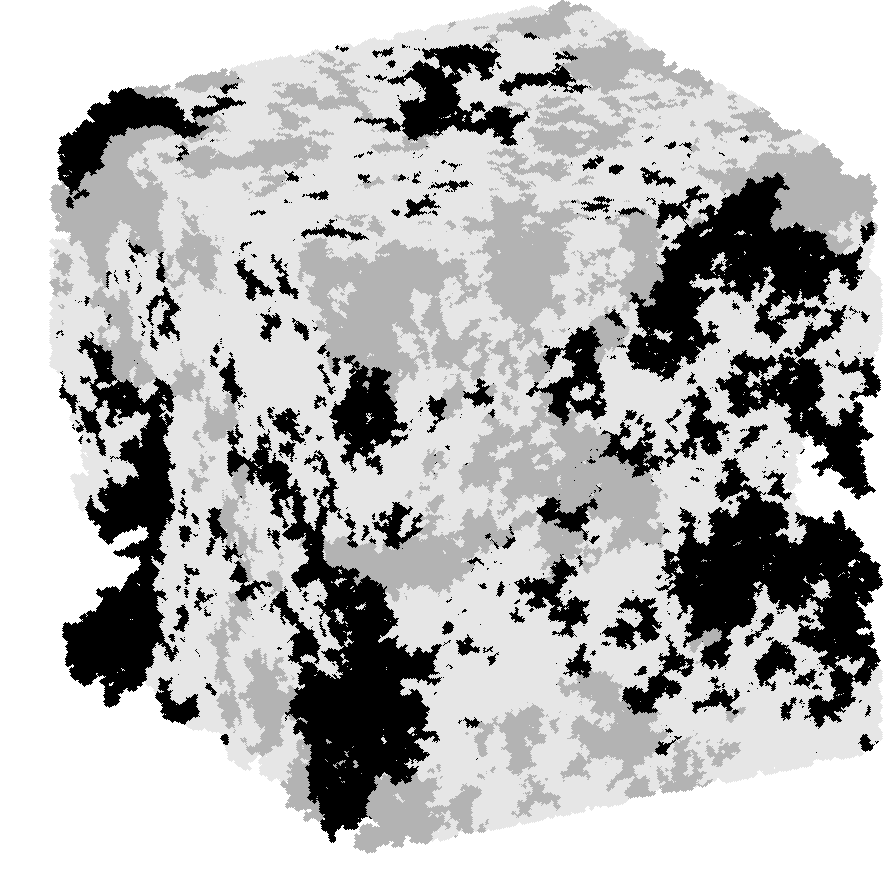}
\caption{3 dimensional PCF \\ {The two images show realisations of PCF on a 256 by 256 by 256 cubic lattice with $\alpha=1$.  Since this is now a 3 dimensional model we have $\alpha<\alpha_c$ and thus we are in the super-critical regime.  The image on the left has the largest warm cluster freezing at $t=0.38$ -- this is shown in black.  The grey cluster is then a  secondary macroscopic cluster which has formed subsequently and wrapped around the first.
The right hand image is the result of the largest warm infinite cluster freezing at $t=0.365$ (black) and the new largest warm cluster freezing at $t=0.4$ (grey).  The subsequent pale grey cluster then leaves us with 3 intertwined macroscopic clusters.
}}
\label{fig:3d}
\end{figure}

Figures \ref{fig:PCF060}, \ref{fig:PCF055} and \ref{fig:PCF050} show the largest clusters from simulations of PCF on a square grid for $\alpha=0.60$, $0.55$ and $0.50$.  These demonstrate the transition from the sub-critical to super-critical regimes.  Figure \ref{fig:3d} then shows simulations of super-critical PCF on a cubic lattice.  Here the freezing times of the largest clusters are such that we get 2 or 3 spanning clusters in the final configuration.  Since the model is not too sensitive to the freezing times of the largest clusters, these images suggest that there is positive probability of 2 of more infinite clusters forming when we run super-critical PCF in dimension 3 or higher.  For further discussion on criticality and the number of infinite clusters the reader is directed to the open problems in Section \ref{sec:problems}.

\paragraph{On a tree}

It turns out that PCF on a tree can be  related to bond percolation on the tree via a suitable time rescaling.  
From this the behaviour of the model on a tree can be well understood, and in particular the distribution of component sizes can be calculated explicitly.

\begin{theorem}\label{thm:cluster_size}
Consider PCF on the rooted binary tree $\mathbf{T}_2$.  There is a critical rate of freezing $\alpha_c=1$ such that for $\alpha<\alpha_c$ infinite clusters form almost surely, but for $\alpha \geq \alpha_c$ all clusters are almost surely finite.  Moreover writing $P_k(\alpha)$ for the probability that the root cluster has size $k$ when it freezes then
\vspace{-0.5\baselineskip}\begin{itemize}
\item if $\alpha<1 \qquad P_k(\alpha)\sim C_{\alpha}\, k^{-2} \qquad$ for some constant $C_{\alpha}$,
\item if $\alpha=1 \qquad P_k(1) \sim C\, k^{-7/4} \qquad$ for some constant $C$,
\item if $\alpha>1 \qquad$then $P_k(\alpha)$ decays exponentially in $k$ at a rate dependent on $\alpha$.
\end{itemize}
\end{theorem}

More generally on a rooted $d$-ary tree we have $\alpha_c=d-1$ and at $\alpha=\alpha_c$ we have $P_k(\alpha_c) \sim C_d\, k^{-(2-1/2d)}$.  The dependence of the exponent on the degree of the graph is rather surprising, and suggests that in contrast to other models there is no lower critical dimension for PCF.

\paragraph{On a complete graph}

Heuristically we can understand PCF on the complete graph $\mathbf{K}_n$ (as $n\longrightarrow\infty$) by comparing with PCF on the $n$-ary tree and rescaling time so that edges open at rate $\frac{1}{n}$.  Our proof then carries over to the mean field case providing a method for rigorously confirming the results of Ben--Naim and Krapivsky, \cite{bn&k1} and \cite{bn&k2}.  See Remark \ref{rem:Knlimitingcase} for more details.

\subsection{Variants of the model}

The PCF process presented here is based upon bond percolation.  It is equally possible to consider a model of site percolation with constant freezing.  In such a model we would open each $v\in\mathcal{V}$ independently at rate $1$, and freeze each open cluster independently at rate $\alpha$.  Once a cluster had frozen all adjacent vertices would remain closed forever.  Theorem \ref{thm:amenable} and Theorem \ref{thm:cluster_size} also hold in this case by following essentially the same proofs.

\subsection{Open problems}\label{sec:problems}

The PCF process leaves us with a wealth of open problems.  Some of these are presented here in the hope of persuading the reader of the richness of the model.

\paragraph{Monotonicity}
Intuitively one would expect that increasing the rate of freezing would lead to clusters freezing more quickly and thus being smaller when they do so.  This intuition might lead us to hypothesise that for $0<\alpha<\beta$ we have
\begin{align}
\mu_{\mathbf{G},\beta} \leq_{\text{st}} \mu_{\mathbf{G},\alpha} \label{eq:stochasticordering}
\end{align}
where $\leq_{\text{st}}$ denotes \emph{stochastic ordering}, and means that $\mu_{\mathbf{G},\beta}(A) \leq \mu_{\mathbf{G},\alpha}(A)$ for all increasing events $A$.  Stochastic ordering is a relatively strong condition and so it could be the case that stochastic ordering does not hold, but that a weaker monotonicity condition does.  However, we do not know of any finite graphs on which (\ref{eq:stochasticordering}) fails.

Such a monotonicity condition is of particular interest because it would imply the existence of a critical value $\alpha_c$.

A related question is whether or not an \emph{FKG inequality} holds for the model.  Whilst a quick check reveals that that the so-called \emph{FKG lattice condition} does not hold even in the case that $\mathbf{G}$ is a line segment with 4 vertices; we have been unable to find a finite graph on which the FKG inequality itself does not hold.  It is worthwhile noting that a BK type inequality does not hold for this model on even the most simple graphs; therefore any attempts to adapt the lace expansion to this model would require more powerful tools. 

\paragraph{Existence of infinite clusters on $\ZZ^d$}

Proposition \ref{prop:no_infinite} tells us that when $\alpha$ is sufficiently large then PCF on $\ZZ^d$ does not produce infinite clusters.  However, simulations (and Figure \ref{fig:PCF050}) suggest that when $\alpha$ is small then an infinite cluster will form.  The following intuition suggests why this should indeed be the case:

Since an edge is joined to at most 2 clusters the probability it is open at time $t$ is at least $p_{\alpha,c}(t)=\frac{1}{1+2\alpha}(1-\ee^{-(1+2\alpha)t})$.  Now let $\alpha>0$ be small, and suppose we run PCF  until some time $t$, where $t$ is sufficiently large for $p_{\alpha,c}(t)>p_c$.  Here $p_c$ is the critical value for bond percolation on $\ZZ^d$.  Thus if the edges were independent then $\ZZ^d$ would contain some infinite cluster.  Now observe that at this time at most $\alpha t$ of the vertices in the graph will be frozen and so -- provided $\alpha$ is sufficiently small -- the number of vertices which have been removed from the infinite percolation cluster (by freezing) are not sufficient to partition the cluster into finite pieces.  Therefore an infinite cluster will remain.

Difficulties in extending this intuition to a proof arise from the fact that the freezing of vertices is not independent, and large clusters can freeze potentially creating barriers which stop warm clusters from growing.  Moreover we should not expect to have good control on the size of these frozen clusters since simulations show that their size has a polynomial tail.  The problem therefore is to find alternative means for showing that an infinite PCF cluster must exist.

\paragraph{Uniqueness of infinite clusters for $d\geq 2$}

During the construction of an infinite volume PCF process it is shown -- in Proposition \ref{prop:unique1} -- that any warm infinite cluster is necessarily unique.  Ben--Naim and Krapivsky showed that this is also true in the mean field case, but that it is also possible for an infinite cluster to freeze and allow subsequent infinite clusters to form.  The geometry of the plane makes it very unlikely that could happen with PCF on $\mathbbm{Z}^2$ -- since any infinite cluster is likely to partition the space into finite pieces.  But there is no such problem in higher dimension, and both simulations (see Figure \ref{fig:3d}) and the following heuristic argument suggest that in dimension $d\geq 3$ it is possible for subsequent infinite clusters to form after previous infinite clusters have frozen.

Let $d\geq 3$ and consider rate $\alpha$ PCF on $\ZZ^d$.  Observe that if we could set $\alpha=0$ then the PCF process would become a percolation process where edges open independently at rate 1.  In this situation it is known that there is some critical time $t_c<\log 2$ at which an infinite cluster will first appear.  Therefore if we have $0<\alpha\ll 1$ it seems reasonable to suppose that the time $t_{\alpha,c}$ at which an infinite cluster appears in the rate $\alpha$ PCF process is close to $t_c$.  Since $\alpha>0$ there is then a positive probability that this infinite cluster will freeze shortly after forming.  Thus provided $\alpha$ is sufficiently small there is a positive probability that an infinite cluster will form and freeze before $t=\log 2$.  At this time the probability that a given edge will be closed and warm is close to $\frac{1}{2}$, so since these edge events will be almost independent there must be an infinite connected set of closed and warm edges in which a new infinite cluster can form.  Thus at some later time $T$ we will have two infinite clusters (one warm and one frozen). 

This leads to the following conjecture.

\begin{conjecture}
Let $d\geq 3$, then provided $\alpha>0$ is sufficiently small there is a positive probability that the final distribution of the rate $\alpha$ PCF process on $\ZZ^d$ will contain 2 or more infinite clusters.
\end{conjecture}

One can show that there can only ever be a finite number of infinite clusters in the final distribution.  However, this  conjecture still leads us to ask how many infinite clusters is it possible to have, and whether or not the maximum number of clusters is controlled by $\alpha$.
Observe that the possibility of having arbitrarily many clusters would imply the possibility of having clusters of arbitrarily low density.  Therefore asking questions about the maximum number of infinite clusters is similar to asking if the percolation function $\theta(p)$ is continuous for higher dimensional lattices.

%

\paragraph{Influence of boundary conditions}
In this paper we use warm boundary conditions to construct an infinite volume limit of the model with free boundary conditions.  We ask therefore what influence the boundary conditions have on the model, and in particular if a infinite volume version of the model where clusters freeze as soon as they touch the boundary exists.

\paragraph{Existence of the model on a general graph}

This paper shows that PCF can be constructed on any countable tree or amenable vertex transitive graph.  We ask therefore if the process can be defined on a more general class of graphs, or indeed on any graph.

%% file: PCF_construction.tex
\section{Constructing PCF}

The construction of PCF relies upon various auxiliary processes and measures.  In order to make the notation as comprehensible as possible a table of notation is included below.

\begin{table}[h]
\centering
    \begin{tabular}{|l|l|l|l|}
    \hline
    Process & Symbol & Measure & Mapping \\ \hline
    PCF process & $(\eta^t)_{t\in[0,\infty]}$  & $\mu=\mu_{\mathbf{G},\alpha}$ & $\psi=\psi_{\GG}:\Pi\longrightarrow\Omega\times[0,\infty]$ \\
    warm PCF process & $(\zeta^t)_{t\in[0,\infty]}$  & $\nu=\nu_{\mathbf{G},\alpha}$ & $\phi=\phi_{\GG}:\Pi\longrightarrow\Omega\times[0,\infty]$ \\
    $\mathrm{liminf}$ of PCF & $(\tilde{\eta}^t)_{t\in[0,\infty]}$  & $\tilde{\nu}=\tilde{\mu}_{\mathbf{G},\alpha}$ & $\tilde{\psi}=\tilde{\psi}_{\GG}:\Pi\longrightarrow\Omega\times[0,\infty]$ \\
    \hline
    \end{tabular}
\end{table}

We now start our construction with a formal definition for the percolation with constant freezing process.

\begin{definition}[PCF on a finite graph]\label{def:PCFonfinite} $ _{}$

Let $\alpha>0$.  Given a finite subgraph $\mathbf{H}=(\mathcal{V}_{\mathbf{H}},\mathcal{E}_{\mathbf{H}})$ of $\GG=(\mathcal{V},\mathcal{E})$, we consider PCF as a Markov process on $\mathbf{H}$ where an edge can be either \emph{open} (state $1$) or \emph{closed} (state $0$), and a vertex can be either \emph{warm} (state $w$) or \emph{frozen} (state $f$).  Thus our state space is $\Omega_{\mathbf{H}}= \{ w,f \}^{\mathcal{V}_{\mathbf{H}}}\times \{0,1\}^{\mathcal{E}_{\mathbf{H}}}$.  Initially each edge is closed, and each vertex is warm.  Therefore $\eta^0 = (w,\ldots,w;0,\ldots,0)$.

Given a configuration $\eta$ we call the subgraph induced by a maximally connected set of open edges a \emph{cluster}.  Note that a path is open irrespective of whether the vertices along it are warm or frozen.  We write $C_v(\eta)=\{w\in\mathcal{V}:w \overset{\eta}{\longleftrightarrow} v \}$ for the cluster containing vertex $v\in\mathcal{V}_{\mathbf{H}}$.  Now for each edge $e\in\mathcal{E}_{\mathbf{H}}$ define $\eta^e$ by
\begin{align}
\eta^e (x) = \left\{ \begin{array}{ll} 1 & x=e \\ \eta(x) & x \neq e , \end{array} \right. \label{eq:^e}
\end{align}
and for each cluster $C\subseteq \mathcal{V}$ define $\eta_{\textsc{c}}$ by
\begin{align}
\eta_{\textsc{c}} (x) = \left\{ \begin{array}{ll} f & x\in C \\ \eta(x) & x \notin C . \end{array} \right. \label{eq:_C}
\end{align}
The generator $G:\Omega_{\mathbf{H}}\times\Omega_{\mathbf{H}}\longrightarrow\mathbbm{R}$ for the PCF process is now defined by
\begin{align}
G(\eta,\eta^e) = 1 \qquad &\text{for all } \eta \text{ and } e=v_1 v_2 \text{ such that } \eta(v_1)=\eta(v_2)=w \nonumber \\
G(\eta,\eta_{\textsc{c}}) = \alpha \qquad &\text{for all clusters } C \text{ of } \eta \label{eq:PCFdef}\\
G(\eta,\theta) = 0 \qquad &\text{otherwise} \nonumber .
\end{align}
\end{definition}

\begin{definition}[PCF measure]\label{def:PCFmeasure}$_{}$

The PCF process described in Definition \ref{def:PCFonfinite} induces a measure $\mu_{\mathbf{H},\alpha}=(\mu_{\mathbf{H},\alpha}^t : t\in[0,\infty))$ on $\Omega_{\mathbf{H}}\times [0,\infty )$.  Here $\mu_{\mathbf{H},\alpha}^t$ is the measure on $\Omega_{\mathbf{H}}= \{ w,f \}^{\mathcal{V}_{\mathbf{H}}}\times \{0,1\}^{\mathcal{E}_{\mathbf{H}}}$ induced by the possible configurations of the PCF process at time $t$.

Since the PCF process on $\mathbf{H}$ can make only finitely many jumps, it will reach some final distribution in finite time.  Therefore as $t\longrightarrow\infty$ the measures $\mu_{\mathbf{H},\alpha}^t$ will converge to $\mu_{\mathbf{H},\alpha}^\infty$ -- the measure of the final distribution.  By including $\mu_{\mathbf{H},\alpha}^\infty$ we extend $\mu_{\mathbf{H},\alpha}$ to the compact space $\Omega_{\mathbf{H}}\times [0,\infty ]$, and call $\mu_{\mathbf{H},\alpha}=(\mu_{\mathbf{H},\alpha}^t : t\in[0,\infty])$ the \emph{PCF measure} of $\mathbf{H}\subseteq\GG$.

In order for us to compare the PCF process on other subgraphs it will be convenient for us to extend the process $(\eta^t)_{t\in[0,\infty]}$ to the whole of $\GG$ by working on the state space $\Omega=\{w,f\}^{\mathcal{V}}\times\{0,1\}^{\mathcal{E}}$ and letting $\eta^t(v)=f$ for all $t$ and each $v\notin\mathcal{V}_{\mathbf{H}}$ and $\eta^t(e)=0$ for all $t$ and each $e\notin\mathcal{E}_{\mathbf{H}}$.  The measure $\mu_{\mathbf{H},\alpha}$ is then extended to $\Omega\times[0,\infty]$ in the obvious way.  Intuitively one might find it helpful to think of this as PCF with free boundary conditions.
\end{definition}

Observe that if we choose an enumeration for the vertices of $\GG$,$\mathcal{V}=\{v_1,v_2,\ldots\}$, then we can identify each cluster by asking which is the vertex of lowest index (or \emph{highest priority}) that it contains.  By doing so we can control the freezing of clusters by attaching independent rate $\alpha$ exponential clocks $X_v$ to each vertex $v\in\mathcal{V}_{\mathbf{H}}$. We can also control the opening of edges with exponential rate $1$ clocks $X_e$ attached to each edge $e\in\mathcal{E}_{\mathbf{H}}$.  This leads to the following algorithm for simulating PCF on a finite graph.

\subsection{An Algorithm for PCF}

\begin{algorithm}[PCF on a finite subgraph]\label{alg:PCF}$_{}$

Let $\alpha>0$, and enumerate the vertices of $\mathcal{V}=\{v_1,v_2,\ldots\}$.  Then given a finite subgraph $\mathbf{H}=(\mathcal{V}_{\mathbf{H}},\mathcal{E}_{\mathbf{H}})\subseteq\GG$ we get $\mathcal{V}_{\mathbf{H}}=\{v_{i_1},\ldots,v_{i_n}\}$.  Now label the vertices of $\mathbf{H}$ by $\ell^0(v_i)=i$ for $1\leq i \leq n$.  These labels will be used to denote which cluster a given vertex is in and satisfy
\begin{align}
\ell^{t}(v_i) = \inf \{ j : v_j\in C_{v_i}(\eta^t) \} = \inf \{ j : \eta^t \text{ contains an open path from } v_j \text{ to } v_i \} \label{eq:labels} .
\end{align}

Let $\{X_v\}_{v\in\mathcal{V}_{\mathbf{H}}}$ be independent $\mathrm{Exp}(\alpha)$ random variables and let $\{X_e\}_{e\in\mathcal{E}_{\mathbf{H}}}$ be independent $\mathrm{Exp}(1)$ random variables.  Our c\`adl\`ag PCF process is now described by its evolution.  At time $t=0$ set $\eta^0=\{w\}^{\mathcal{V}} \times \{0\}^{\mathcal{E}}$, and for each $x\in\mathcal{V}\cup\mathcal{E}$ at time $t=X_x$
\vspace{-0.5\baselineskip}\begin{enumerate}
\item[(a)] if $t=X_{v_i}$ then set $C=\{v_k\in\mathcal{V}:\ell^{t^-}(v_k)=i\}$ and let $\eta^{t}=(\eta^{t^-})_{\textsc{c}}$.
\item[(b)] if $t=X_e$ with $e=v_i v_j \in \mathcal{E}$, then if $\eta^{t^-}(v_i)=\eta^{t^-}(v_j)=1$ set $\eta^{t}=(\eta^{t^-})^e$.

To ensure that (\ref{eq:labels}) continues to hold we set $\ell^{t}(v_k)=\mathrm{min}\{\ell^{t^-}(v_i),\ell^{t^-}(v_j)\}$ for all $v_k$ with $\ell^{t^-}(v_k)\in\{\ell^{t^-}(v_i),\ell^{t^-}(v_j)\}$.
\end{enumerate}
\end{algorithm}

This algorithm can be thought of as a map $\psi_{\mathbf{H}}$ from $\{X_v\}_{v\in\mathcal{V}_{\mathbf{H}}}\cup\{X_e\}_{e\in\mathcal{E}_{\mathbf{H}}}$ to a PCF process $(\eta^t)_{t\in[0,\infty]}$.  From this it is clear that a PCF event $A\subseteq\Omega_{\mathbf{H}}$ is measurable with respect to $\mathcal{F}_{\mathbf{H}}=\sigma(\{X_v\}_{v\in\mathcal{V}_{\mathbf{H}}}\cup\{X_e\}_{e \in\mathcal{E}_{\mathbf{H}}})$.  Note that here an explicit choice of enumeration of the vertices is required.  However, since the $X_v$ are independent it is clear that this choice is arbitrary and does not affect the law of the process itself.  Therefore it is clear that this algorithm gives a process satisfying (\ref{eq:PCFdef}), and so as $\mathbf{H}$ is finite we see that this gives the process of Definition \ref{def:PCFonfinite}.

\subsection{The warm PCF process}

Given an amenable vertex transitive graph $\GG$ we begin the construction of the PCF process on $\GG$ by first constructing an auxiliary process which obeys a suitable monotonicity condition.  We consider a finite subgraph $\mathbf{H}\subseteq\GG$ for now and define \emph{warm PCF} as follows.

\begin{definition}[Warm PCF on a finite subgraph]\label{def:warmPCF} $ _{}$

Suppose $\mathbf{H}=(\mathcal{V}_{\mathbf{H}},\mathcal{E}_{\mathbf{H}})$ is a finite subgraph of $\GG=(\mathcal{V},\mathcal{E})$, we say that $v\in \mathcal{V}_{\mathbf{H}}$ is a \emph{boundary vertex} if 
it meets is an edge $e=vw$ in $\mathcal{E}_{\mathbf{G}}$ which is not in $\mathcal{E}_{\mathbf{H}}$.
We denote the set of boundary vertices by $\partial\mathbf{H}$.  Now given a configuration $\zeta\in\Omega_{\mathbf{H}} = \{w,f\}^{\mathcal{V}_{\mathbf{H}}} \times \{0,1\}^{\mathcal{E}_{\mathbf{H}}}$ we say that $C$ is a \emph{boundary cluster} of $\zeta$ if $C$ is a cluster of $\zeta$ with $C\cap\partial\mathbf{H}\neq\emptyset$. We say that a cluster $C\subseteq\mathbf{H}$ of $\zeta$ which is not a boundary cluster is an \emph{interior cluster}.

We now set our initial configuration to $\zeta^0 = \{w\}^{\mathcal{V}_{\mathbf{H}}} \times \{0\}^{\mathcal{E}_{\mathbf{H}}}$, and define the generator $G$ for the warm PCF process on $\mathbf{H}$ to be
\begin{align}
G(\zeta,\zeta^e) = 1 \qquad &\text{for all }\zeta \text{ and } e=v_1 v_2 \text{ such that } \zeta(v_1)=\zeta(v_2)=w \nonumber\\
G(\zeta,\zeta_{\textsc{c}}) = \alpha \qquad &\text{for all interior clusters } C \text{ of } \zeta \label{eq:warmPCFdef} \\
G(\zeta,\theta) = 0 \qquad &\text{otherwise}. \nonumber
\end{align}
Therefore boundary clusters will remain warm forever.  Intuitively one might like to think of this as PCF with warm wired boundary conditions.

It will be useful for us to be able to compare warm PCF on two subgraphs $\mathbf{H}_1$ and $\mathbf{H}_2$.  To make sense of this we extend warm PCF on $\mathbf{H}$ to a process on $\GG$ by working with the state space $\Omega=\{ w,f \}^{\mathcal{V}}\times\{0,1\}^{\mathcal{E}} $ and setting our initial configuration to $\zeta^0 = \{w\}^{\mathcal{V}} \times \{0\}^{\mathcal{E}_{\mathbf{H}}} \times \{1\}^{\mathcal{E} \setminus \mathcal{E}_{\mathbf{H}}}$.  Since $G$ only acts on a finite number of edges and vertices we see that (\ref{eq:warmPCFdef}) gives a well defined process.
\end{definition}

Warm PCF differs from PCF in that the transition $\zeta\longrightarrow\zeta_{\textsc{c}}$ only occurs for interior clusters in (\ref{eq:warmPCFdef}), while in PCF this transition can take place with all clusters (see (\ref{eq:PCFdef})).  This means that the warm PCF process dominates the PCF process -- a notion we shall make rigorous in Lemma \ref{lem:dominate}.

\begin{definition}[Warm PCF measure]\label{def:warmPCFmeasure}$_{}$

In the same way as described in Definition \ref{def:PCFmeasure} the warm PCF process induces a measure $\nu_{\mathbf{H},\alpha}^{}=(\nu_{\mathbf{H},\alpha}^t : t\in[0,\infty])$ on $\Omega\times [0,\infty ]$.  This is the \emph{warm PCF measure} of $\mathbf{H}\subseteq\GG$.
\end{definition}

We also remark that Algorithm \ref{alg:PCF} can be adapted to warm PCF by replacing \textit{(a)} with
\vspace{-0.5\baselineskip}\textit{\begin{enumerate}
\item[(a')] if $t=X_{v_i}$ then set $C=\{v_k\in\mathcal{V}:\ell(v_k)=i\}$ and provided $C\cap\partial\mathbf{H}=\emptyset$ let $\zeta^{t}=(\zeta^{t^-})_{\textsc{c}}$.
\end{enumerate}}
We shall denote this modified algorithm as Algorithm \ref{alg:PCF}*.

Again this algorithm can be thought of as giving a function $\phi_{\mathbf{H}}$ from $\{X_v\}_{v\in\mathcal{V}_{\mathbf{H}}}\cup\{X_e\}_{e\in\mathcal{E}_{\mathbf{H}}}$ to a warm PCF process $(\zeta^t)_{t\in[0,\infty]}$.  This then tells us that  a warm PCF event $A\subseteq\Omega$ is also measurable with respect to the $\mathcal{F}=\sigma(\{X_v\}_{v\in\mathcal{V}}\cup \{ X_e \}_{e\in\mathcal{E}})$.

\begin{definition}[Coupling]\label{def:coupling}$_{}$

Given a countable graph $\GG=(\mathcal{V},\mathcal{E})$ and a collection of finite subgraphs $\{\mathbf{H}_i\}_{i\in I}$ we can use Algorithm \ref{alg:PCF} and Algorithm \ref{alg:PCF}* to couple the PCF and warm PCF processes on each of the $\mathbf{H}_i$.  To do this we fix an enumeration of the vertices of $\GG$, and assign an $\mathrm{Exp}(\alpha)$ random variable $X_v$ to each $v\in \mathcal{V}$ and an $\mathrm{Exp}(1)$ random variable $X_e$ to each $e\in \mathcal{E}$.  Denote the probability space for the random variables $\{X_x\}_{x\in\GG}$ by $(\Pi,\lambda_{\mathbf{G},\alpha})$.  Now for each $i\in I$ Algorithm \ref{alg:PCF} and Algorithm \ref{alg:PCF}* give us coupled mappings $\psi_i:(\Pi,\lambda_{\mathbf{G},\alpha})\longrightarrow(\Omega\times[0,\infty],\mu_{\mathbf{H}_i,\alpha})$ and $\phi_i:(\Pi,\lambda_{\mathbf{G},\alpha})\longrightarrow(\Omega\times[0,\infty],\nu_{\mathbf{H}_i,\alpha})$ from $\{X_x\}_{x\in\GG}$ to a PCF process on $\mathbf{H}_i$ and a warm PCF process on $\mathbf{H}_i$.
\end{definition}

The key reason for introducing the intermediate warm PCF process is that it obeys a monotonicity condition.  Given two configurations $\eta,\zeta\in\Omega$ we say $\eta\leq\zeta$ if and only if $\zeta(e)=1$ for all $e\in\mathcal{E}$ with $\eta(e)=1$ and $\zeta(v)=w$ for all $v\in\mathcal{V}$ with $\eta(v)=w$.  That is, $\zeta$ has more open edges and warmer vertices than $\eta$. From this partial ordering of $\Omega$ we define $A\subseteq\Omega$ to be an \emph{increasing event} if and only if $\eta\in A$ and $\eta\leq\zeta$ implies $\zeta\in A$.  Whence we get a \emph{stochastic ordering} by saying that if $\mu$ and $\nu$ are two measures on $\Omega$ then $\mu\leq_{\text{st}}\nu$ if and only if $\mu(A)\leq\nu(A)$ for all measurable increasing events $A$.

To extend this notion of stochastic ordering to processes on $\GG$ we consider trajectories $(\eta^{t})_{t\in [0,\infty]}$ and $(\eta^{t})_{t\in [0,\infty]}$, and say that $(\eta^{t})\leq(\zeta^{t})$ if and only if $\eta^{t}\leq\zeta^{t}$ for all $t\in[0,\infty]$.  An event $A\subseteq\Omega\times [0,\infty]$ is then an \emph{increasing process event} if and only if $(\eta^{t})\in A$ and $(\eta^{t})\leq(\zeta^{t})$ implies $(\zeta^{t})\in A$.  Therefore given two process $P$ and $Q$  with respective measures $\mu=(\mu^t)_{t\in[0,\infty]}$ and $\nu=(\nu^t)_{t\in[0,\infty]}$, we say that $Q$ \emph{stochastically dominates} $P$ if and only if $\mu(A)\leq\nu(A)$ for all measurable increasing process events $A$, and write $\mu\leq_{\text{st}}\nu$.

\begin{lemma}[Monotonicity for warm PCF]\label{lem:monotonicity}
Suppose $\mathbf{H}_1=(\mathcal{V}_1,\mathcal{E}_1)$ and $\mathbf{H}_2=(\mathcal{V}_2,\mathcal{E}_2)$ are finite subgraphs of $\GG$ with $\mathbf{H}_1 \subseteq \mathbf{H}_2$, then the rate $\alpha$ warm PCF process on $\mathbf{H}_1$ stochastically dominates the rate $\alpha$ warm PCF process on $\mathbf{H}_2$.
\bp
Here we work with the processes extended to $\GG=(\mathcal{V},\mathcal{E})$, and use the notion of coupling from Definition \ref{def:coupling}.  Let $\{X_v\}_{v\in\mathcal{V}}$ and $\{X_e\}_{e\in\mathcal{E}}$ be as in Definition \ref{def:coupling} and let $(\zeta_1^{t})_{t\in [0,\infty]}$ and $(\zeta_2^{t})_{t\in [0,\infty]}$ be the trajectories associated with  warm PCF on $\mathbf{H}_1$ and $\mathbf{H}_2$ respectively.  It now suffices to show that the coupling gives $\zeta_1^{t} \geq \zeta_2^{t}$ for all $t\in[0,\infty]$.

For the purposes of our proof we shall show further that if $C$ is a warm interior cluster of $\zeta_1^{t}$ then it is also a warm interior cluster of $\zeta_2^{t}$.  Since $\mathbf{H}_1 \subseteq \mathbf{H}_2$ this is clear at $t=0$, and because $\mathcal{E}_1 \subseteq\mathcal{E}_2$ we also have $\zeta_1^0 = \{w\}^{\mathcal{V}} \times \{0\}^{\mathcal{E}_{1}} \times \{1\}^{\mathcal{E} \setminus \mathcal{E}_{1}} \geq \{w\}^{\mathcal{V}} \times \{0\}^{\mathcal{E}_{2}} \times \{1\}^{\mathcal{E} \setminus \mathcal{E}_{2}} = \zeta_2^0$.

Now because $(\zeta_1^{t})$ and $(\zeta_2^{t})$ can only change at discrete times $t=X_x$, where $x\in\mathcal{V}_2 \cup \mathcal{E}_2$, it suffices for us to show that if our hypothesis holds and we have $\zeta^{t^-}_1\geq\zeta_2^{t^-}$ at time $t^-$ then the hypothesis also holds at time $t$ for $t\in\{X_x : x\in\mathcal{V}_2\cup\mathcal{E}_2\}$.  There are several cases:
\vspace{-0.5\baselineskip}\begin{enumerate}
\item[(a)] If $t=X_{v}$ and $v$ is the vertex of highest priority in a warm interior cluster $C$ of $\zeta_1^{t^-}$ then $C$ must also be a warm interior cluster of $\zeta_2^{t^-}$ with $v$ again being the vertex of highest priority.  So we have $\zeta_1^{t}=(\zeta_1^{t^-})_{\textsc{c}}\geq (\zeta_2^{t^-})_{\textsc{c}}=\zeta_2^{t}$, and $C$ is no longer a warm interior cluster of $\zeta_1^{t}$ or $\zeta_2^{t}$.

If $v$ is a vertex of highest priority in a warm interior cluster $C$ of $\zeta_2^{t^-}$ but no int $\zeta_1^{t^-}$ then we get $\zeta_1^{t}=\zeta_1^{t^-}\geq\zeta_2^{t^-} \geq (\zeta_2^{t^-})_{\textsc{c}} = \zeta_2^{t}$.  Moreover, $C$ remains warm in $\zeta_1^t$ but not in $\zeta_2^t$.

If $v$ is not the vertex of highest priority in some warm interior cluster of either $\zeta_1^{t^-}$ or $\zeta_2^{t^-}$ then $\zeta_1^{t}$ and $\zeta_2^{t}$ remain unchanged.
\item[(b)] If $t=X_e$ with $e=v_i v_j$, then if $\zeta_2^{t^-}(v_i)=\zeta_2^{t^-}(v_j)=w$ we must also have $\zeta_1^{t^-}(v_i)=\zeta_1^{t^-}(v_j)=w$, since $\zeta_1^{t^-} \geq\zeta_2^{t^-}$ and so any vertex which is warm in $\zeta_2^{t^-}$ must also be warm in $\zeta_1^{t^-}$.  Therefore if $e$ is open in $\zeta_2^t$ then it must also be open in $\zeta_1^t$, and so we see that $\zeta_1^{t} \geq \zeta_2^{t}$.

If $C_{v_i}=C_{v_j}$ then there is no change in the clusters so we assume they are disjoint.  If $C_{v_i}$ and $C_{v_j}$ are warm interior clusters of $\zeta_1$ then they are also warm interior clusters of $\zeta_2$ and so $C_{v_i}\cup C_{v_j}$ is a warm interior cluster of both $\zeta_1$ and $\zeta_2$.  If one of $C_{v_i}$ or $C_{v_j}$ is not a warm interior clusters of $\zeta_1^{t^-}$, then it must be a must be a boundary cluster of $\zeta_1^{t^-}$ and so $C_{v_i}$ and $C_{v_j}$ will be (contained in) boundary clusters of $\zeta_1^{t}$, and so our condition on the clusters is satisfied.  

If $\zeta_2^{t^-}(v_i)=\zeta_2^{t^-}(v_j)=w$ but $\zeta_2^{t^-}(v_i)=f$ or $\zeta_2^{t^-}(v_j)=f$ then we have $\zeta_1^{t}=(\zeta_1^{t^-})^e\geq\zeta_1^{t^-} \geq \zeta_2^{t^-} = \zeta_2^{t}$.  Moreover, one of $C_{v_i}$ or $C_{v_j}$ must have been a boundary cluster of $\zeta_1^{t^-}$ and so $C_{v_i}$ and $C_{v_j}$ will be (contained in) boundary clusters of $\zeta_1^{t}$, and so again the condition on the clusters is satisfied. 

If $\zeta_1^{t^-}(v_i)=f$ or $\zeta_1^{t^-}(v_j)=f$ then $\zeta_2^{t^-}(v_i)=f$ or $\zeta_2^{t^-}(v_j)=f$ and so $\zeta_1^{t}$ and $\zeta_2^{t}$ remain unchanged, and the condition remains satisfied.
\qedhere \end{enumerate}
\ep
\end{lemma}

Note that no such monotonicity condition need exist  when we compare PCF on finite subgraphs $\mathbf{H}_1\subseteq\mathbf{H}_2$.  However, we do obtain monotonicity when we compare the warm PCF process to the PCF process on a finite subgraph $\mathbf{H}\subseteq\GG$.

\begin{lemma}[Stochastic domination of PCF by warm PCF]\label{lem:dominate}
Suppose $\mathbf{H}\subseteq\GG$ is a finite subgraph, then the rate $\alpha$ warm PCF process on $\mathbf{H}$ stochastically dominates the rate $\alpha$ PCF process on $\mathbf{H}$.  Moreover if $(\eta^{t})_{t\in[0,\infty]}=\psi((s_x)_{x\in\GG})$ and $(\zeta^{t})_{t\in[0,\infty]}=\phi((s_x)_{x\in\GG})$ are coupled PCF and warm PCF processes, then at all $t\in[0,\infty]$ any interior cluster $C$ of $\zeta^t$ is also a cluster of $\zeta^t$.
\bp
This is proved in the same way as Lemma \ref{lem:monotonicity} by using the coupling of Definition \ref{def:coupling}.  The details are omitted.
\ep
\end{lemma}

Given $\alpha>0$ we can now construct the unique infinite volume warm PCF process on an amenable vertex transitive graph $\GG$.  Let $\GG_1\subseteq\GG_2\subseteq\ldots$ be any exhaustion of $\GG$ with each $\GG_n$ finite, and for each $n$ consider the measure $\nu_n$ induced by warm PCF on $\GG_n$.  By monotonicity we have that $\nu_1,\nu_2,\ldots$ converges to some measure $\nu$.

To check the limit $\nu$ doesn't depend on the exhaustion suppose $\mathbf{H}_1\subseteq\mathbf{H}_2\subseteq\ldots$ is another exhaustion of $\GG$ with each $\mathbf{H}_n$ finite, and let $\nu^{*}_n$ be the induced warm PCF measure for each $n$.  Again the $\nu^{*}_n$ must converge to some measure $\nu^{*}$.  Now since $\{\GG_n:n\in\mathbbm{N}\}$ and $\{\mathbf{H}_n:n\in\mathbbm{N}\}$ are both exhaustions of $\GG$ there exists subsequences $(i_k)_{k\geq 1}$ and $(j_k)_{k\geq 1}$ such that $\GG_{i_1}\subseteq\mathbf{H}_{j_1}\subseteq\GG_{i_2}\subseteq\mathbf{H}_{j_2}\subseteq\ldots$, and from this we see that $\nu_{i_1}\geq_{\text{st}}\nu^{*}_{j_1}\geq_{\text{st}}\nu_{i_2} \geq_{\text{st}}\nu^{*}_{j_2}\geq_{\text{st}}\ldots$ and deduce that this sequence converges.  Since this sequence contains a subsequence of both $(\nu_n)$ and $(\nu^{*}_n)$ we must have that $\nu=\nu^{*}$.  Therefore we define this unique limiting $\nu=\nu_{\GG,\alpha}$ to be the \emph{infinite volume warm PCF measure} on $\GG$.

Because of the way $\nu$ has been constructed by taking an exhaustion $\GG_1\subseteq\GG_2\subseteq\ldots$ of $\GG$ it is clear that for any finite $\La\subseteq\GG$ the restriction of warm PCF on $\GG_n$ to $\La$ will converge to a unique limiting process as $n\longrightarrow\infty$.  This limiting process is the restriction of$\nu$ to $\La$.  Therefore $\nu$ defines an infinite volume warm PCF measure in the sense of local limits.

Suppose $\{X_v\}_{v\in\mathcal{V}}$ and $\{X_e\}_{e\in\mathcal{E}}$ are independent exponential random variables at respective rates $\alpha$ and $1$, and denote the product space of these by $\Pi$ with measure $\lambda=\lambda_{\GG,\alpha}$.  In Definition \ref{def:coupling} we introduced functions $\phi_n:(\Pi,\lambda)\longrightarrow(\Omega\times[0,\infty],\nu_n)$ given by Algorithm \ref{alg:PCF}* for each finite $\GG_n$.  We have now seen that the $\phi_n$ have an almost sure limit as $n\longrightarrow\infty$ which we define to be $\phi_{\GG}=\phi:(\Pi,\lambda)\longrightarrow(\Omega\times[0,\infty],\nu)$.

\begin{lemma}[Translation invariance]\label{lem:translation}
Suppose that $f:\GG\longrightarrow\GG$ is an isometry, and that $A$ is measurable with respect to the $\sigma$-algebra of $\nu$, then $\nu(A)= \nu(f(A))$.  Thus if $\GG$ is vertex transitive then $\nu$ is translation invariant.
\bp
Consider a fixed exhaustion of $\GG$, $\GG_1\subseteq\GG_2\subseteq\ldots$ say, and let  $\nu_1,\nu_2,\ldots$ be the induced warm PCF measures.  Now given any isometry $f$ of $\GG$ let $\nu^{*}_1,\nu^{*}_2,\dots$ be the measures induced by warm PCF on $f(\GG_1),f(\GG_2),\ldots$.  For any measurable event $A$ we have that $\nu_n(A)=\nu^{*}_n(f(A))$ for all $n$.  Therefore we get
\begin{align}\label{eq:translation}
\nu(A)=\lim_{n\rightarrow\infty} \nu_{n}(A) = \lim_{n\rightarrow\infty} \nu^{*}_n(f(A)) = \nu(f(A)) ,
\end{align}
as required.
\ep
\end{lemma}

Observe that up to now we have only required for $\GG$ to be countable -- in order for the exhaustion $\GG_1\subseteq\GG_2\subseteq\ldots$ to exist; and for $\GG$ to be locally finite -- in order for us to make sense of boundary vertices.  Thus we have proved the following.

\begin{proposition}
For every locally finite countable graph $\GG$ and every fixed rate of freezing $\alpha>0$, there exists an infinite volume warm PCF process on $\GG$ in the sense of local limits.  Moreover in the case that $\GG$ is vertex transitive then $\nu$ is also translation invariant. 
\end{proposition}

\subsection{Proof of Theorem \ref{thm:amenable}}

Whilst this auxiliary  process is interesting in its own right, it differs from the infinite volume PCF process of Theorem \ref{thm:amenable} in that an infinite cluster will never freeze.  This is because an infinite cluster must have been a boundary cluster of each $\GG_n$.  Our strategy therefore is to attach an exponential clock to each warm infinite cluster -- making it freeze at rate $\alpha$.  To do this we must first show that there are only finitely many warm infinite clusters at any given time $T$.  This is where it is necessary to use the fact that $\GG$ is amenable.  The content of Proposition \ref{prop:unique1} below is to show that any warm infinite cluster is necessarily unique.  Once we know this we can complete our construction by freezing the unique warm infinite cluster (if it exists) at times $0<T_1<T_2<\ldots$.  Here the $T_i$ will be chosen to have independent rate $\alpha$ exponential increments.  The way we freeze the warm infinite cluster at time $T_i$ is defined below.

\begin{definition}[Freezing the infinite cluster at times $\mathcal{T} = \{0<T_1<T_2<\ldots\}$]\label{def:freezeinfinite}$_{}$

Let $0<T_1<T_2<\ldots$ be a fixed sequence of times, and let $\GG_1\subseteq\GG_2\subseteq\ldots$ be an exhaustion $\GG$ with each $\GG_n$ finite.  Now use the random variables $\{X_e\}_{e\in\mathcal{E}}$ and $\{X_v\}_{v\in\mathcal{V}}$ to run coupled warm PCF on each $\GG_n$ -- giving a processes $(\zeta^{t}_n)_{t\in [0,\infty]}$ -- and let $(\zeta^{t}_{})_{t\in [0,\infty]}$ be the limiting warm PCF process on $\GG$.

At time $T_1$ set $C_{\infty,1}= \{ v\in\mathcal{V} : v \text{ is contained in an infinite cluster of } \zeta^{T_1^-}\text{ and } \zeta^{T_1^-}(v)=w \}$, and for each $n$ let $C_{n,\infty,1} = C_{\infty,1} \cap \GG_n$.  Note that if no infinite cluster exists at time $T_1$ then $C_{\infty,1} = \emptyset$.  We now change each warm PCF process to have $\zeta^{T_1}_n = (\zeta^{T_1^-}_n)_{C_{n,\infty,1}}$, before continuing to use Algorithm \ref{alg:PCF}* to evolve $\zeta^{t}_{n}$ for $t>T_1$.  Denote these new processes by $(\zeta_{n,T_1}^{t})_{t\in [0,\infty ]}$.

Suppose $m\leq n$, then since $\GG_m\subseteq \GG_n$ we have $C_{m,\infty,1}\subseteq C_{n,\infty,1}$, and so because $\zeta^{T_1^-}_{n}\leq\zeta^{T_1^-}_m$ we must also have $\zeta^{T_1}_{n}= (\zeta^{T_1^-}_n)_{C_{n,\infty,1}} \leq (\zeta^{T_1^-}_m)_{C_{m,\infty,1}} = \zeta^{T_1}_m$.  Moreover, since freezing $C_{n,\infty,1}$ and $C_{m,\infty,1}$ does not change the warm interior clusters of $\zeta^{T_1}_n$ or $\zeta^{T_1}_m$ the proof of Lemma \ref{lem:monotonicity} tells us that $\zeta_{m,T_1}^{t} \geq \zeta_{n,T_1}^{t}$ for all $t$.

Because this monotonicity property still holds we can repeat the argument given above to see that $(\zeta_{n,T_1}^{t})_{t\in [0,\infty ]} \longrightarrow (\zeta_{T_1}^{t})_{t\in [0,\infty ]}$ as $n\longrightarrow\infty$.  This in turn gives us a new measure $\nu_{T_1}$.

Now we repeat this modification at times $T_2,T_3,\ldots$ by setting $C_{\infty,j+1}= \{ v\in\mathcal{V} : v$ is contained in an infinite cluster of $\zeta_{T_1,\ldots,T_{j}}^{T_{j+1}^-}$ and $\zeta_{T_1,\ldots,T_{j}}^{T_{j+1}^-}(v)=w \}$ and $C_{n,\infty,j+1} = C_{\infty,j+1} \cap \GG_n$ for each $n$.  Doing so we obtain a sequence of new measures $\nu_{T_1,T_2}$, $\nu_{T_1,T_2,T_3}$, etc. Since $\nu_{T_1,\ldots,T_m}$ and $\nu_{T_1,\ldots,T_n}$ are equal when restricted to $t\in [0,T_m]$, if we make the further assumption that $T_j \longrightarrow\infty$ as $j\longrightarrow\infty$ then these measures must converge to some final measure $\nu_{\mathcal{T}}$.

Note that for each $\mathcal{T}=\{0<T_1<T_2<\ldots\}$ with $T_j\longrightarrow\infty$ we can combine Algorithm \ref{alg:PCF}* with the procedure above to obtain a function $\phi_{\mathcal{T}}:(\Pi,\lambda)\longrightarrow(\Omega\times[0,\infty],\nu_{\mathcal{T}})$ taking a sequence $(s_x)_{x\in\GG}\in\Pi$ to a process $(\zeta_{\mathcal{T}}^{t})_{t\in[0,\infty]}\in\Omega\times[0,\infty]$.  
\end{definition}

\begin{remark}
If $\GG$ is vertex transitive then for fixed $\mathcal{T}=\{0<T_1<T_2,\ldots\}$ with $T_j\longrightarrow\infty$ the measure $\nu_{\mathcal{T}}$ is again translation invariant.  This can be seen inductively since the set $C_{\infty,1}$ is determined by the translation invariant measure $\nu^{T_1}$, and so by applying the argument of Lemma \ref{lem:translation} to $\GG\setminus C_{\infty,1}$ we see that $\nu_{T_1}$ is translation invariant for $t>T_1$ (as well as for $0\leq t \leq T_1$).  Now suppose we know that $\nu_{T_1,\ldots,T_j}$ is translation invariant, then because $C_{\infty,j+1}$ is determined by $\nu_{T_1,\ldots,T_j}^{T_j+1}$ it must be a translation invariant set, and so by applying Lemma \ref{lem:translation} to $\GG\setminus C_{\infty,1}\cup\ldots\cup C_{\infty,j+1}$ we see that $\nu_{T_1,\ldots,T_{j+1}}$ is again a translation invariant measure.
\end{remark}

Suppose multiple disjoint warm infinite clusters were present at $T_j^-$, then the process of Definition \ref{def:freezeinfinite} would lead to them all freezing at $T_j$ and so their freezing would not be independent.  However, we shall now show that any warm infinite cluster is unique almost surely, and thus the independence of freezing is maintained.

\begin{proposition}[Uniqueness of the warm infinite clusters]\label{prop:unique1}$_{}$

Let $\alpha>0$, and suppose $\GG$ is an amenable vertex transitive graph.  Consider the infinite volume warm PCF process on $\GG$ at rate $\alpha$, then at any time $T\geq 0$ we have that $\nu^{T}$--almost every $\zeta\in\Omega$ contains at most one infinite cluster.

Moreover, if we modify the process so that any warm infinite cluster freezes at fixed times $0<T_1<\ldots<T_k$, as in Definition \ref{def:freezeinfinite}, then at any time $T\geq T_k$ we have that $\nu_{T_1,\ldots,T_k}^{T}$--almost every warm infinite cluster is unique.

\begin{remark}
A theorem of Burton and Keane, \cite[Theorem 2]{burtonkeane}, says that for any measure $\mu$ on $\Omega$ which has the so-called \emph{finite energy property} a configuration $\zeta$ has at most one infinite cluster $\mu$--almost surely.  Whilst it is possible to show that $\nu^{T}$ meets the hypotheses of \cite[Theorem 2]{burtonkeane}, the modified measures $\nu_{T_1,\ldots,T_k}^{T}$ do not (changing the state of one vertex from warm to frozen can have infinite effect).  The proof of this proposition deals with this problem by using only events for which a finite energy condition does hold.
\end{remark}

\bp
Observe that the first claim is a special case of the second with $k=0$, so we set $k\geq 0$ and let $\mathcal{T}=\{0=T_0<T_1<\ldots<T_k\}$ be fixed.  The claim is trivial if $T=T_k$, so we also fix $T > T_k$.  Given a configuration $\zeta$ define $N(\zeta)\in\{0,1,2,\ldots\}\cup\{\infty\}$ to be the number of warm infinite clusters of $\zeta$.  If $\nu_{\mathcal{T}}^{T}$ was ergodic then we would know that $N$ is equal to some constant $\tilde{N}$ almost surely.  Our plan of attack would then be to show that
\vspace{-0.5\baselineskip}\begin{enumerate}
\item If $\tilde{N}\in\{2,3,\ldots\}$ then $\nu_{\mathcal{T}}^{T}(\mathbf{A}_1)>0$, where $\mathbf{A}_1$ is an event whose positive probability would contradict $N=\tilde{N}$ $\nu_{\mathcal{T}}^{T}$--almost surely.
\item If $\tilde{N}=\infty$ then $\nu_{\mathcal{T}}^{T}(\mathbf{A}_2)>0$, where $\mathbf{A}_2$ is another event whose positive probability would again lead to a contradiction.
\end{enumerate}
We would then be able to conclude that $N(\eta)\in\{0,1\}$ $\nu_{\mathcal{T}}^{T}$--almost surely.

However, since we do not know that $\nu_{\mathcal{T}}^{T}$ is ergodic, we have to rely on the Ergodic Decomposition Theorem.  This says that a translation invariant measure (such as $\nu_{\mathcal{T}}^{T}$) can be decomposed into ergodic measures.  More precisely, for any translation invariant measure space $\Omega$ there exists a measurable map from $\Omega$ to the space of ergodic measures on $\Omega$, $m:\Omega\longrightarrow\mathscr{E}(\Omega)$, such that
\begin{align}
\mu(A)=\int_{x\in\Omega}m_x (A)\,\dd \mu(x) ,
\end{align}
for all translation invariant measures $\mu$ and all measurable $A$.  For details see \cite[Proposition 10.26]{kallenberg} or \cite{ergodic}.  Moreover, the following lemma allows us to relate the changes in energy of $m_x$ and $\nu_{\mathcal{T}}^{T}$ from the modification of a finite set of edges and vertices.

\begin{lemma}\label{lem:decomposition}
Suppose $\La\subseteq\GG$ is a finite set of edges and vertices, and $A\in\Omega$ is an event depending only on $\La$.  Write $\mathcal{F}_{\mathbf{G}\setminus\La}$ for the $\sigma$-algebra generated by the edges and vertices of  $\mathbf{G}\setminus\La$.  Then for $\nu_{\mathcal{T}}^{T}$--almost all $m_x$ in the ergodic decomposition of $\nu_{\mathcal{T}}^{T}$ and each $A\in\mathcal{F}_{\La}$ we have that 
\begin{align*}
m_x(A \,|\, \mathcal{F}_{\mathbf{G}\setminus\La}) = \nu_{\mathcal{T}}^{T}(A \,|\, \mathcal{F}_{\mathbf{G}\setminus\La}) \qquad m_x \text{--almost everywhere} .
\end{align*}
\vspace{-\baselineskip}\bp
The result follows from a careful adaptation of the proof of \cite[Lemma 1]{gkn}.
%
%
\ep
\end{lemma}

\begin{remark}
One might be tempted to try and overcome the lack of ergodicity by working with the product space $(\Pi,\lambda)$ rather than $(\Omega,\nu_{\mathcal{T}})$.  However, whilst the event $\{$there are $k$ warm infinite clusters at time $T\}$ is clearly tail measurable with respect to $\Omega$, we do not know that its preimage $\phi^{-1}_{\mathcal{T}}(\{$there are $k$ warm infinite clusters at time $T\})$ is tail measurable with respect to $\Pi$.  Therefore whilst we shall use $(\Pi,\lambda)$ later on for making explicit calculations -- such as in (\ref{eq:Step1}) -- it is necessary for us begin the proof by working with $(\Omega,\nu_{\mathcal{T}})$.
\end{remark}

\paragraph*{Step 1}
For contradiction suppose that $\nu_{\mathcal{T}}^{T}(N(\zeta)\in\{2,3,\ldots\})>0$, then there must be some $k\in\{2,3,\ldots\}$ with $\nu_{\mathcal{T}}^{T}(N(\zeta)=k)>0$.  For each $v\in\mathcal{V}$ let $d(0,v)$ be the length of the shortest path from $0$ to $v$, and set $\La_r = \{ v\in\mathcal{V}  :  d(0,v) \leq r \}$.  We now define
\begin{align}
A_r = \{ \zeta : N(\zeta)=k \text{ and }\zeta \text{ contains at least 2 warm infinite clusters intersecting }\La_r\}.
\end{align}
Now $\nu_{\mathcal{T}}^{T}(A_r)\longrightarrow \nu_{\mathcal{T}}^{T}(N(\zeta)=k)>0$ as $r\longrightarrow\infty$, therefore we can fix $R\in\mathbbm{N}$ with $\nu_{\mathcal{T}}^{T}(A_R)>0$.  For a given configuration $\zeta$ let
\begin{align*}
\La_R^{*} (\zeta) &= \{v\in\mathcal{V} : v \text{ is in } \La_R \text{ or in a finite cluster of } \zeta \text{ intersecting } \La_R \} \\
\bar{\La}_R(\zeta) &= \La_R^{*} \cup \{ e\in\mathcal{E} : e \text{ meets }\La_R^{*}(\zeta) \} .
\end{align*}
Observe that $\nu_{\mathcal{T}}^{T}(A_R \cap \{ \bar{\La}_R (\zeta) \subseteq \La_{\rho}\})\longrightarrow \nu_{\mathcal{T}}^{T}(A_R)$ as $\rho\longrightarrow\infty$, and therefore $\nu_{\mathcal{T}}^{T}(A_R \cap \{ \bar{\La}_R (\zeta) \subseteq \La_{\rho}\})>0$ for some $\rho\in\mathbbm{N}$.  Moreover, because $\La_{\rho}$ is finite there must be some subset $\La\subseteq\La_{\rho}$ with $\nu_{\mathcal{T}}^{T}(A_R \cap \{ \bar{\La}_R (\zeta) = \La\})>0$.  These sets are shown in Figure \ref{fig:maketrifurcation}.  We now fix such a $\La$ and write $A_{R,\La} = A_R \cap \{ \bar{\La}_R (\zeta) = \La\}$.

Because $\nu_{\mathcal{T}}^{T}(A_{R,\La})>0$ there must be a positively measurable set of $x\in\Omega$ such that $m_x$ is in the ergodic decomposition of $\nu_{\mathcal{T}}^{T}$ with $m_x(A_{R,\La})>0$.  We fix such an $m_x$, and note that since $m_x$ is ergodic and translation invariant it must have $N(\zeta)=k$, $m_x$--almost surely.  To get our contradiction we shall let $\mathbf{A}_1$ be the event that $N(\zeta)<k$; we claim that
\begin{align}
m_x(\mathbf{A}_1)>0. \label{eq:claim1}
\end{align}

\bpc
To estimate $m_x(\mathbf{A}_1)$ we define two more events:
\begin{align}
B_{R,\La} &= \{ \tilde{\zeta}\in\Omega : \text{there is some } \zeta\in A_{R,\La} \text{ with } \tilde{\zeta}(x)=\zeta(x) \text{ for all } x\in \GG\setminus\La\} , \label{eq:Br} \\
O_{R,\La} &= \{\zeta\in\Omega : \text{all vertices in }\La\text{ are warm, all edges with two ends in }\La \text{ are open}\nonumber\\
&\qquad\text{and all edge with only one end in }\La\text{ are closed} \} \label{eq:Or} .
\end{align}
Note that $B_{R,\La}\in\mathcal{F}_{\GG\setminus\La}$ and $O_{R,\La}\in\mathcal{F}_{\La}$.

Suppose $\tilde{\zeta}\in B_{R,\La}\cap O_{R,\La}$, then we can find $\zeta\in A_{R,\La}$ with $\tilde{\zeta}(x)=\zeta(x)$ for all $x\in\GG\setminus\La$.  Now $\zeta$ has $k$ warm infinite clusters, at least two of which meet $\La_R$.  Therefore since all the edges of $\La_R$ are open and warm, $\tilde{\zeta}$ will have at least two of the infinite clusters of $\zeta$ joined together.  Moreover, since $\zeta$ and $\tilde{\zeta}$ are equal outside $\La$ no new infinite clusters are created.  Thus $N(\tilde{\zeta})<k$, and so $B_{R,\La}\cap O_{R,\La}\subseteq\mathbf{A}_1$.  Because $B_{R,\La}\in\mathcal{F}_{\GG\setminus\La}$ we can apply Lemma \ref{lem:decomposition} to see that
\begin{align}
m_x(\mathbf{A}_1)\geq m_x(B_{R,\La}\cap O_{R,\La}) &= \ex (m_x(1_B \, 1_O | \mathcal{F}_{\GG\setminus\La})) \nonumber \\
&= \ex(1_B \,m_x(1_O | \mathcal{F}_{\GG\setminus\La})) \nonumber \\
&= \ex(1_B \,\nu_{\mathcal{T}}^{T}(1_O | \mathcal{F}_{\GG\setminus\La})) .
\end{align}
where $B=B_{R,\La}$ and $O=O_{R,\La}$.  Since $m_x(B_{R,\La})>0$ we can prove the claim by showing that $\nu_{\mathcal{T}}^{T}(O_{R,\La}|\zeta_{\GG\setminus\La})>0$ for each feasible configuration $\zeta_{\GG\setminus\La}$ of the edges and vertices of $\GG\setminus\La$.

Let $\zeta_{\GG\setminus\La}$ be fixed.  We now switch to working with the product space $(\Pi,\lambda)$.  Using the function defined in Definition \ref{def:freezeinfinite}, $\phi_{\mathcal{T}} : (\Pi,\lambda)\longrightarrow(\Omega\times [0,\infty], \nu_{\mathcal{T}}^{T})$, we have
\begin{align*}
\nu_{\mathcal{T}}^{T}(O_{R,\La}|\zeta_{\GG\setminus\La}) = \lambda (\phi_{\mathcal{T}}^{-1}(O_{R,\La}) | \phi_{\mathcal{T}}^{-1} (\zeta_{\GG\setminus\La}) ) ,
\end{align*}
where $\phi_{\mathcal{T}}^{-1}(\zeta_{\GG\setminus\La}) = \{ (s_x)_{x\in\GG} : \phi_{\mathcal{T}} ((s_x)_{x\in\GG})^T(x) = \zeta_{\GG\setminus\La}(x) \text{ for all } x\in\GG\setminus\La\}$.  From this we define $U_{\zeta_{\GG\setminus\La}}$ to be the set of $(t_x)_{x\in\GG}$ for which there exists $(s_x)_{x\in\GG}\in \phi_{\mathcal{T}}^{-1}(\zeta_{\GG\setminus\La})$ with $t_x=s_x$ for all $x\in\GG\setminus\La$, and note that $U_{\zeta_{\GG\setminus\La}}\supseteq\phi_{\mathcal{T}}^{-1}(\zeta_{\GG\setminus\La})$.  Without loss of generality we may assume that the vertices of $\La$ have highest priority in the algorithm used to construct $\phi_{\mathcal{T}} $.  Now define
\begin{align}
U_{\zeta_{\GG\setminus\La}}^O &= U_{\zeta_{\GG\setminus\La}} \cap \{(s_x) : s_v > T \text{ for all } v\in\La \} \cap \{ (s_x) : s_e > T \text{ for all } e \text{ with one end in } \La \} \nonumber \\
&\qquad \cap\{ (s_x) : s_e < T \text { for all } e \text{ with two ends in } \La \} .
\end{align}
Note that each of the events on the right hand side are independent of $U_{\zeta_{\GG\setminus\La}}$.  We now include two lemmas which allow us to analyse $\phi_{\mathcal{T}}((s_x)_{x\in\GG})$ for $(s_x)_{x\in\GG}\in U_{\zeta_{\GG\setminus\La}}$.

\begin{lemma}\label{lem:conditionaledges}
Set $T>0$ and let $E\subseteq\mathcal{E}$ and $V\subseteq\mathcal{V}$.  Suppose that $(s_x)_{x\in\GG}\in\Pi$ is such that setting $\phi((s_x)_{x\in\GG})=\zeta$ gives $\zeta^{T}(e)=0$ for all $e\in E$ and $\zeta^{T}(v)=w$ for all $v\in V$.  If $(\tilde{s}_x)_{x\in\GG}$ satisfies $\tilde{s}_x=s_x$ for all $x\in\GG\setminus (E\cup V)$, and $\tilde{s}_x>T$ for all $x\in E\cup V$, then writing $\phi((\tilde{s}_x)_{x\in\GG})=\tilde{\zeta}$ gives $\tilde{\zeta}^{T}=\zeta^{T}$ (and indeed $\tilde{\zeta}^t=\zeta^t$ for all $t\in [0,T]$).

Moreover, if $\mathcal{T}$ is fixed then the result still holds when $\phi$ is replaced by $\phi_{\mathcal{T}}$.
\bp
Since $\phi$ (respectively $\phi_{\mathcal{T}}$) is a limit of the $\phi_n^{}$ (respectively $\phi_{n,\mathcal{T}}$) it suffices to consider the algorithm for warm PCF on one of the $\GG_n$.  Suppose $e\in\mathcal{E}_n$ with $s_e < T$ then since $\zeta^{T}(e)=0$ the transition at time $s_e$ is $\zeta^{s_e^-}\longrightarrow\zeta^{s_e}=\zeta^{s_e^-}$.  Similarly if $v\in\mathcal{V}_n$ with $s_v < T$ then since $\zeta^{T}(v)=w$ the transition at time $s_v$ is also $\zeta^{s_v^-}\longrightarrow\zeta^{s_v}=\zeta^{s_v^-}$. Therefore since $\tilde{\zeta}$ does not change at $s_x$ we have that $\tilde{\zeta}^{s_x^-}=\zeta^{s_x^-}$ implies $\tilde{\zeta}^{s_x}=\zeta^{s_x}$, and so the lemma is proved by induction.
\ep
\end{lemma}

\begin{lemma}\label{lem:conditionalvertices}
Set $T>0$, let $V\subseteq\mathcal{V}$, and suppose our enumeration of $\mathcal{V}$ is such that if $v_i\in V$ and $v_j\in\mathcal{V}\setminus V$ then $i<j$ -- i.e. $V$ has higher priority than $\mathcal{V}\setminus V$.  Now given $(s_x)_{x\in\GG}\in\Pi$ let $(\tilde{s}_x)_{x\in\GG}$ satisfy $s_x=\tilde{s}_x$ for all $x\in\GG\setminus V$ and $\tilde{s}_x>T$ for all $x\in V$.  If we write $\phi((s_x)_{x\in\GG})=\zeta$ and $\phi((\tilde{s}_x)_{x\in\GG})=\tilde{\zeta}$ then for all $t\in [0,T]$ we  have
\vspace{-0.5\baselineskip}\begin{enumerate}
\item[(a)] $\tilde{\zeta}^t \geq \zeta^t$,
\item[(b)] for each $e\in\mathcal{E}$ with $\zeta^t(e)=0$ and $\tilde{\zeta}^t(e)=1$ there exists an open path in $\tilde{\zeta}^t$ from $e$ to $V$.
\end{enumerate}
Moreover, if $\mathcal{T}=\{0<T_1<T_2<\ldots\}$ is fixed then the result still holds if we replace $\phi$ by $\phi_{\mathcal{T}}$.
\bp
We consider our exhaustion $\GG_1\subseteq \GG_2\subseteq \ldots$ of $\GG$.  Since $\phi$ is a limit of the $\phi_n$, to show the result it suffices to show that the result holds for each $\phi_n$.  Now consider the algorithm for warm PCF on $\GG_n$.  If a cluster $C$ intersects $V$ then its freezing time is controlled by a vertex $v\in V$, therefore since $s_v>T$ for all $v\in V$ we can treat each $v\in V$ as a boundary vertex up until time $T$.  \textit{(a)} then follows by using the coupling argument of Lemma \ref{lem:monotonicity} with $\mathbf{H}_1=\GG_n\setminus V$ and $\mathbf{H}_2=\GG_n$.

Suppose an edge $e\in\mathcal{E}$ has $\zeta^t(e)=0$ and $\tilde{\zeta}^t(e)=1$ for some $t\in [0,T]$, then since $\zeta^t(e)=0$ we must have $e\in\GG_n$, and as  $\tilde{\zeta}^t(e)=1$ we must have $s_e=\tilde{s}_e<t$.  Now there can only be finitely many such edges (as $\GG_n$ is finite), $\{ e_1,\ldots,e_k \}$ say, and so to prove \textit{(b)} we can consider what must happen at each of the times $s_{e_i}$ in order.

At a time $t=\min\{s_{e_1},\ldots,s_{e_k}\}$ the edges of $\zeta^{t^-}$ and $\tilde{\zeta}^{t^-}$ are equal and thus so are the clusters, but then $\zeta$ has the transition $\tilde{\zeta}^{t^-}\longrightarrow\tilde{\zeta}^{t}=(\tilde{\zeta}^{t^-})^e$ whilst $\zeta$ is unchanged.  Thus one of the end points of $e$ must have $\zeta^{t}(v)=f$ and $\tilde{\zeta}^{t}(v)=w$.  The only way this can happen is if the freezing of $v$ is controlled by a vertex in $V$ and so there must be a path from $e$ to $V$.  At subsequent times $t\in\{s_{e_1},\ldots,s_{e_k} \}$ we must also have some vertex $v$ with $\zeta^{t}(v)=f$ and $\tilde{\zeta}^{t}(v)=w$.  This can only happen if either the freezing of $v$ is controlled by some vertex in $V$ or if the freezing of $v$ was controlled by different vertices in $\zeta$ and $\tilde{\zeta}$.  The first case implies that there is an open path from $e$ to $V$, and the second case implies that $C_v$ must differ in $\zeta$ and $\tilde{\zeta}$.  For $C_v$ to differ, $v$ must be joined to some edge $e'$ with $\zeta^t(e')=0$ and $\tilde{\zeta}^t(e')=1$.  Therefore $e'\in\{e_1,\ldots,e_k\}$ with $s_{e'}<s_e$, and so we can apply induction to see that there is an open path from $e$ to $V$ (through $e'$).

The result for $\phi_{\mathcal{T}}$ follows by the same argument.
\ep
\end{lemma}

Suppose $(t_x)_{x\in\GG}\in U_{\zeta_{\GG\setminus\La}}^O$ and find $(s_x)_{x\in\GG}\in\phi_{\mathcal{T}}^{-1}(\zeta_{\GG\setminus\La})$ with $t_x=s_x$ for all $x\in\GG\setminus\La$.  Write $\zeta=\phi_{\mathcal{T}}((s_x)_{x\in\GG})$ and $\tilde{\zeta} = \phi_{\mathcal{T}}((t_x)_{x\in\GG})$, and define $E=\{ e\in\La : e\text{ has exactly one end in }\La\}$ and $V=\{v\in\partial\La_R : v \text{ is in a warm infinite cluster of }\zeta_{\GG\setminus\La}\}$.  Suppose $(\tau_x)_{x\in\GG}$ has $\tau_x=t_x$ for all $x\in E\cup V$ and $\tau_x=s_x$ for all $x\in\GG\setminus E\cup V$, then by applying Lemma \ref{lem:conditionaledges} we see that $\phi_{\mathcal{T}}((\tau_x)_{x\in\GG})=\zeta$.  Now since the edges on the boundary of $\La$ are closed in $\phi_{\mathcal{T}}((\tau_x)_{x\in\GG})$ we can use Lemma \ref{lem:conditionalvertices} with $(\tau_x)_{x\in\GG}$ and $(t_x)_{x\in\GG}$ to deduce that $\tilde{\zeta}(x)=\zeta(x)$ for all $x\in\GG\setminus\La$.  
A simple check reveals that $\tilde{\zeta}^T \in O_{R,\La}$, therefore $U^O_{\zeta_{\GG\setminus\La}}\subseteq \phi_{\mathcal{T}}^{-1}(O_{R,\La})\cap U_{\zeta_{\GG\setminus\La}}$.  We can now estimate $\nu_{\mathcal{T}}^{T}(O_{R,\La} |\zeta_{\GG\setminus\La}) $,
\begin{align}
\nu_{\mathcal{T}}^{T}(O_{R,\La} |\zeta_{\GG\setminus\La}) &= \lambda(\phi_{\mathcal{T}}^{-1}(O_{R,\La})|\phi_{\mathcal{T}}^{-1}(\zeta_{\GG\setminus\La})) \nonumber \\
&\geq \lambda(U_{\zeta_{\GG\setminus\La}}^O | U_{\zeta_{\GG\setminus\La}}) \nonumber \\
&= (\ee^{-\alpha T})^{ |\La\cap\mathcal{V}|} \times (\ee^{-T})^{|\La\cap\mathcal{E}\setminus E|} \times (1-\ee^{-T})^{|E|} > 0 \label{eq:Step1}.
\end{align}
Thus we indeed have $m_x(\mathbf{A}_1)>0$, and so (\ref{eq:claim1}) is proved.
\epc

\paragraph*{Step 2}
We now move on to showing that $N(\zeta)\neq\infty$ $\nu_{\mathcal{T}}^{T}$--almost surely.  To do this we recall the notion of \emph{trifurcation points}.  Given a configuration $\zeta$ we say a vertex $v\in\mathcal{V}$ is a trifurcation point if $v$ is contained in some infinite component $C$ and removing $v$ (and its adjoining edges) from $C$ leaves three disjoint infinite clusters $C_1$, $C_2$ and $C_3$ with $C_1 \cup C_2 \cup C_3 \cup \{ v\} = C$.  By following the argument of Burton and Keane \cite{burtonkeane} we see that any translation invariant measure $\mu$ on an amenable vertex transitive graph $\GG$ can not have $\mu (0 \text{ is a trifurcation point}) > 0$.

\begin{figure}[p]
\centering
\includegraphics[width=0.6\textwidth]{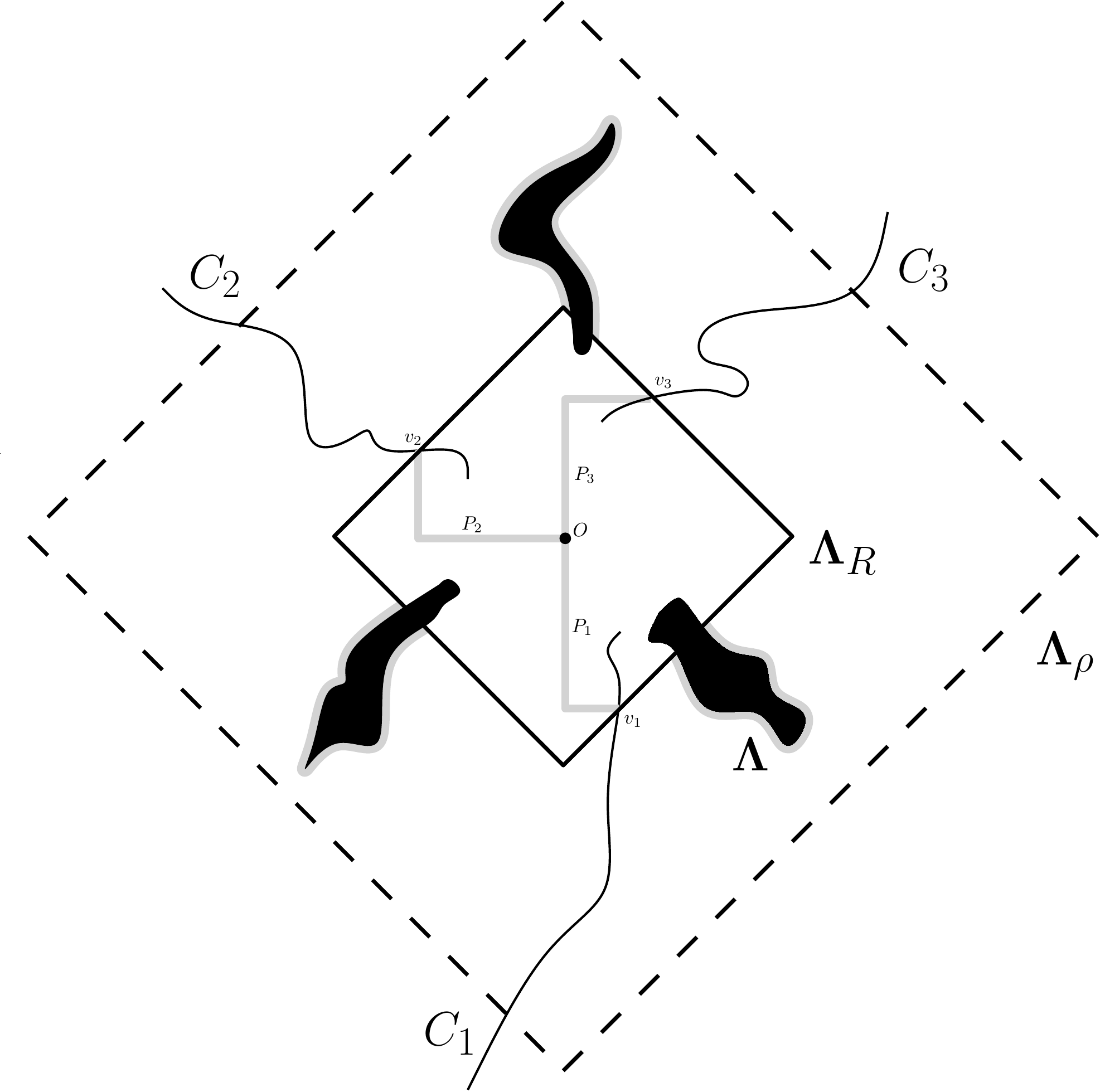}
\caption{The diagram shows a box $\La_R$ intersecting three infinite clusters.  $\La$ is formed with the addition of the finite clusters (and their boundary edges) which meet $\La_R$.  Note that $\rho$ has been chosen so that $\La\subseteq\La_{\rho}$ with positive probability. \\
Also shown the three paths $P_1,P_2,P_3$ from $v_1,v_2,v_3$ to $O$.  In the second part of the proof we ensure that $P=P_1 \cup P_2 \cup P_3$ is open so that the origin becomes a trifurcation point.}
\label{fig:maketrifurcation}
\end{figure}

\begin{figure}[p]
\centering
\includegraphics[width=0.6\textwidth]{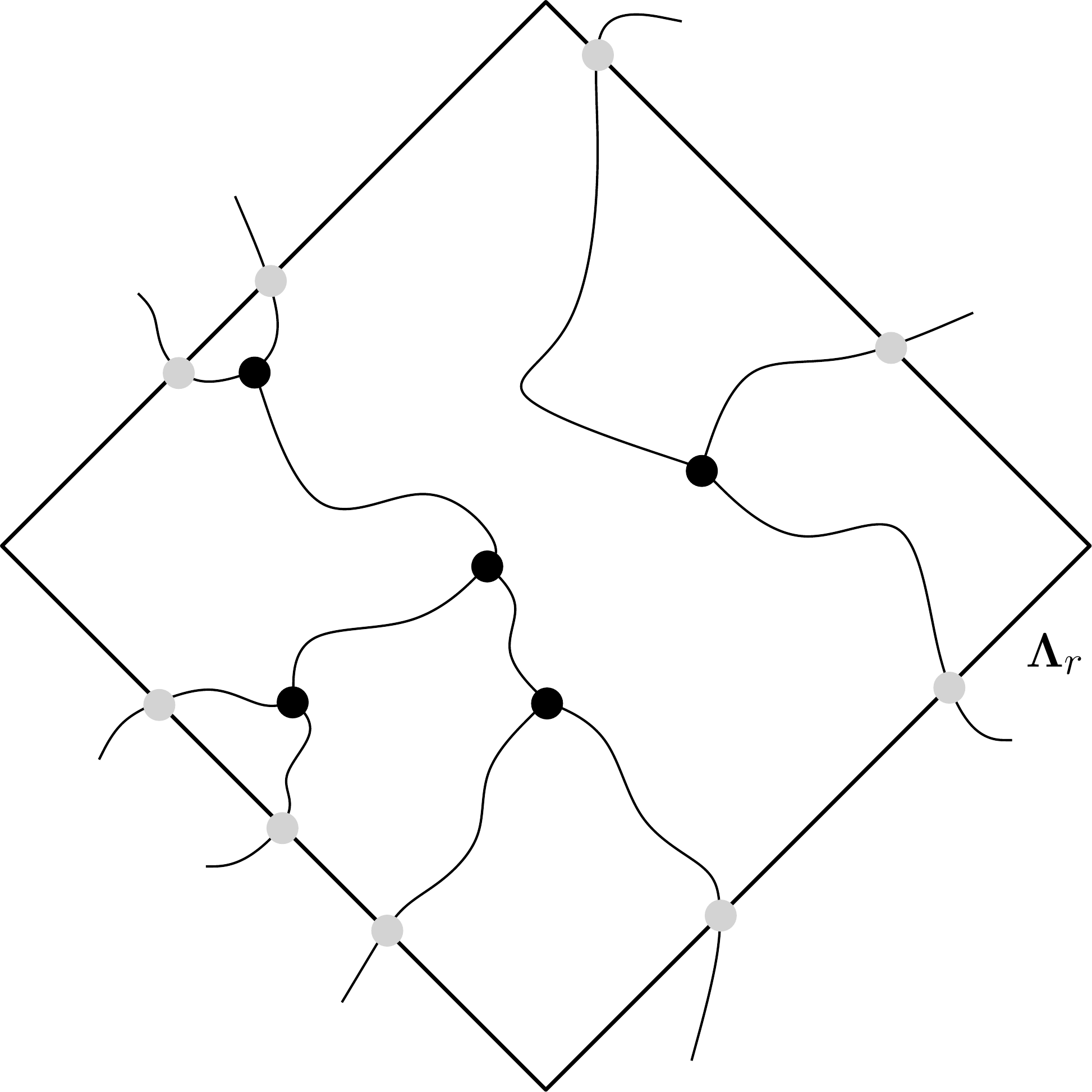}
\caption{This digram shows $\La_r$ containing five trifurcation points in two infinite clusters.  Observe that any infinite cluster containing $n$ trifurcation points must intersect the boundary of $\La_r$ in at least $n+2$ places.}
\label{fig:trifurcation}
\end{figure}

The argument of \cite{burtonkeane} (which is repeated in \cite{grimmett1}) is based on the observation that if there are $n$ trifurcation points in $\La_{r}$ then $\partial\La_r = \La_r \setminus \La_{r-1}$ must intersect infinite clusters in at least $n+2$ places.  This is illustrated in Figure \ref{fig:trifurcation}.  The Ergodic Theorem then tells us that if a translation invariant measure $\mu$ satisfies $\mu (0 \text{ is a trifurcation point}) > 0$ then there exists $c>0$ such that 
\begin{align} \mu (\La_r \text{ contains at least } c | \La_r | \text{ trifurcation points}  ) \longrightarrow 1 \quad \text{as} \quad r \longrightarrow \infty \label{eq:amenable}.
\end{align}
Because $\GG$ is amenable then $|\partial\La_r | / |\La_r |\longrightarrow 0$ as $r\longrightarrow\infty$.  This means that $|\partial\La_r|<c|\La_r|$ for sufficiently large $r$, and so (\ref{eq:amenable}) is impossible.  Therefore we let $\mathbf{A}_2$ be the event that $0$ is a trifurcation point, and get a contradiction by showing that it occurs with positive probability.

\begin{claim'}
If $\nu_{\mathcal{T}}(N(\zeta)=\infty)>0$ then $\nu_{\mathcal{T}}(\mathbf{A}_2)>0$.
\bpc
We prove the claim in a similar way to before:  Find a finite $\La$ where a suitable configuration on $\GG\setminus\La$ occurs with positive probability, and modify the edges of $\La$ with only a finite cost.  However, unlike before no ergodicity assumption is required.  We start by defining $A_r^{\infty}$ in a similar way to $A_r$, but this time insisting that $\La_r$ meets at least three warm infinite clusters.
\begin{align}
A_r^{\infty} = \{ \zeta : N(\zeta)=\infty \text{ and }\zeta \text{ contains at least 3 warm infinite clusters intersecting }\La_r\}.
\end{align}
Since $\nu_{\mathcal{T}}(N(\zeta)=\infty)>0$ then there exists some $R\in\mathbbm{N}$ with $\nu_{\mathcal{T}}(A^{\infty}_R)>0$.  As before we can now find a finite $\La\supseteq\La_R$ for which $A^{\infty}_{R,\La}=A^{\infty}_R\cap \{\bar{\La}_R(\zeta)=\La\}$ has $\nu_{\mathcal{T}}(A^{\infty}_{R,\La})>0$.  Given three distinct vertices $v_1,v_2,v_3\in\partial\La_R$ we define a further event
\begin{align}
A^{\infty}_{R,\La,v_1,v_2,v_3} = A^{\infty}_{R,\La} \cap \{ \zeta : v_1,v_2,v_3 \text{ are in disjoint warm infinite clusters of } \zeta \} .
\end{align}
There are only finitely many distinct 3-sets $\{v_1,v_2,v_3\}$, and so with positive probability there exists $v_1,v_2,v_3 \in\partial\La_R$ with $\nu_{\mathcal{T}}(A^{\infty}_{R,\La,v_1,v_2,v_3})>0$.  Fix such a 3-set and let  $P_1,P_2,P_3\subseteq\La_R$ be disjoint paths from $v_1,v_2,v_3$ to 0, with $P=P_1\cup P_2\cup P_3$.  See Figure \ref{fig:maketrifurcation}.  We now define analogues to $B_{R,\La}$ and $O_{R,\La}$.
\begin{align}
B_{R,\La,v_1,v_2,v_3} &= \{ \tilde{\zeta}\in\Omega : \text{there is some } \zeta\in A_{R,\La,v_1,v_2,v_3} \text{ with } \tilde{\zeta}(x)=\zeta(x) \text{ for all } x\in \GG\setminus\La\} ,  \\
\Delta_{R,\La,v_1,v_2,v_3} &= \{\zeta \in\Omega: \text{all vertices in }\La\text{ are warm, all edges } e\in\La\setminus P \text{ are closed } \nonumber \\
&\qquad \text{and all edges }e\in P \text{ are open} \} .
\end{align}
Note that since $A^{\infty}_{R,\La,v_1,v_2,v_3}\subseteq B_{R,\La,v_1,v_2,v_3}$ we have $\nu_{\mathcal{T}}(B_{R,\La,v_1,v_2,v_3})>0$. Moreover, $B_{R,\La,v_1,v_2,v_3}$ depends only on the edges and vertices of $\GG\setminus\La$ and $\Delta_{R,\La,v_1,v_2,v_3}$ depends only on the edges and vertices of $\La$.

Suppose $\tilde{\zeta}\in B_{R,\La,v_1,v_2,v_3}\cap \Delta_{R,\La,v_1,v_2,v_3}$, then we can find $\zeta$ in $A^{\infty}_{R,\La,v_1,v_2,v_3}$ with $\zeta(x)=\tilde{\zeta}(x)$ for all $x\in \GG\setminus\La$.  Now since $\zeta$ contains three disjoint infinite clusters $C_1,C_2,C_3$ meeting $v_1,v_2,v_3$, and because $v_1,v_2,v_3$ are joined to $0$ by open paths in $\tilde{\zeta}$, then $\tilde{\zeta}$ must contain an infinite cluster $C$ with a trifurcation point at $0$.  This means $\mathbf{A}_2 \subseteq B_{R,\La,v_1,v_2,v_3}\cap \Delta_{R,\La,v_1,v_2,v_3}$, and so
\begin{align*}
\nu_{\mathcal{T}}(\mathbf{A}_2) \geq \nu_{\mathcal{T}}(\Delta_{R,\La,v_1,v_2,v_3}| B_{R,\La,v_1,v_2,v_3}) \times \nu_{\mathcal{T}}(B_{R,\La,v_1,v_2,v_3}) .
\end{align*}
Therefore to prove the claim it suffices to show that $\nu_{\mathcal{T}}(\Delta_{R,\La,v_1,v_2,v_3}| \zeta_{\GG\setminus\La})>0$ for each feasible configuration $\zeta_{\GG\setminus\La}$.  To do this we again switch to working with the product space $(\Pi,\lambda)$.  Using our measure preserving function $\phi_{\mathcal{T}} : (\Pi,\lambda)\longrightarrow(\Omega\times [0,\infty], \nu_{\mathcal{T}}^{T})$ we have
\begin{align*}
\nu_{\mathcal{T}}^{T}(\Delta_{R,\La,v_1,v_2,v_3}|\zeta_{\GG\setminus\La}) = \lambda (\phi_{\mathcal{T}}^{-1}(\Delta_{R,\La,v_1,v_2,v_3}) | \phi_{\mathcal{T}}^{-1} (\zeta_{\GG\setminus\La}) ) .
\end{align*}
Define $U_{\zeta_{\GG\setminus\La}}$ to be the set of $(t_x)_{x\in\GG}$ for which there exists $(s_x)_{x\in\GG}\in\phi_{\mathcal{T}}^{-1}(\zeta)$ with $t_x = s_x$ for all $x\in\GG\setminus\La$, and let
\begin{align*}
U_{\zeta_{\GG\setminus\La}}^{\Delta} &= U_{\zeta_{\GG\setminus\La}}\cap \{ (s_x) : s_v>T \text{ for all } v\in\La \} \cap \{(s_x):s_e<T \text{ for all }e\text{ in }P\} \\
&\qquad \cap \{(s_x):s_e > T \text{ for all } e\in\La\setminus P \} .
\end{align*}
Using Lemma \ref{lem:conditionaledges} and Lemma \ref{lem:conditionalvertices} we can check that $U_{\zeta_{\GG\setminus\La}}^{\Delta}\subseteq \phi_{\mathcal{T}}^{-1}(\Delta_{R,\La,v_1,v_2,v_3})\cap U_{\zeta_{\GG\setminus\La}}$, and so we can estimate $\nu_{\mathcal{T}}(\Delta_{R,\La,v_1,v_2,v_3}|\zeta_{\GG\setminus\La})$ as follows.
\begin{align}
\nu_{\mathcal{T}}(\Delta_{R,\La,v_1,v_2,v_3}|\zeta_{\GG\setminus\La}) &= \lambda(\phi_{\mathcal{T}}^{-1}(\Delta_{R,\La,v_1,v_2,v_3})|\phi_{\mathcal{T}}^{-1}(\zeta_{\GG\setminus\La})) \nonumber \\
&\geq \lambda(U_{\zeta_{\GG\setminus\La}}^{\Delta} | U_{\zeta_{\GG\setminus\La}}) \nonumber \\
&= (\ee^{-\alpha T})^{ |\La\cap\mathcal{V}|} \times (\ee^{-T})^{|P|} \times (1-\ee^{-T})^{|\La\cap\mathcal{E}\setminus P|} > 0.
\end{align}
Thus we can conclude that $\nu_{\mathcal{T}}(\mathbf{A}_2)>0$.
\epc
\end{claim'}

As both $\nu_{\mathcal{T}}(N(\zeta)\in\{2,3,\ldots\})>0$ and $\nu_{\mathcal{T}}(N(\zeta)=\infty)>0$ lead to a contradiction, we must have $N(\zeta)\in\{0,1\}$ $\nu_{\mathcal{T}}$--almost surely.
\ep
\end{proposition}

In Definition \ref{def:freezeinfinite} we freeze each of the infinite clusters present some time $T_i$.  However, since we have shown that a warm infinite cluster is necessarily unique, then freezing all warm infinite clusters at $T_i$ will result in at most one infinite cluster freezing.  Therefore we let $X_{\infty,1},X_{\infty,2},\ldots$ be a sequence of independent $\mathrm{Exp}(\alpha)$ random variables, set $T_i=\sum_{j=1}^{i} X_{\infty,i}$ for $i=1,2,\ldots$, denote $\mathcal{T}=\{ T_1,T_2,\ldots\}$, and then define $\mu = \nu_{\mathcal{T}} $.

It is easy to see that the process corresponding to $\mu$ satisfies the transition rates (\ref{eq:PCFdef}) of Definition \ref{def:PCFonfinite}.  Indeed each warm edge $e$ will be contained in $\GG_n$ for $n$ sufficiently large, and so will open at rate $1$.  Likewise each warm cluster $C$ is either finite and thus an interior cluster of some $\GG_n$ meaning it freezes at rate $\alpha$ or it is infinite meaning its freezing is controlled by some $X_{\infty,i}$ -- again at rate $\alpha$.  We shall now conclude the proof of Theorem \ref{thm:amenable} by showing that the measure $\mu$ corresponds to PCF on $\GG$ in the sense of local limits.

\begin{theorem}
Fix $\alpha>0$ and let $\mathbf{G}$ be an amenable vertex transitive graph.  Then if $\La\subseteq\GG_1\subseteq\GG_2\subseteq ...$ (each $\GG_n$ finite) is an exhaustion of $\GG$ with PCF measures $\mu_{1},\mu_{2},\ldots$ then $\mu_n|_{\La}\longrightarrow\mu|_{\La}$ as $n\longrightarrow\infty$.  Thus $\mu_n\longrightarrow\mu$ in the sense of local limits.  Here $.|_{\La}$ denotes the restriction of a measure to the subgraph $\La$ and $\mu = \nu_{\mathcal{T}}$ is as defined above.
\bp
Since $\Omega$ is a product of discrete spaces and thus compact the obvious approach is to use Prokhorov's Theorem to extract a convergent subsequence $(\mu_{n_k})_{k\geq 1}$, let $\mu^* = \lim_{k\rightarrow\infty} \mu_{n_k}$ and compare $\mu^*|_{\La}$ to $\mu|_{\La}$.  However, since $\mu^*$ depends on the subsequence $(n_k)_{k\geq 1}$ the is no a priori reason why this limit should be unique.  We get around this problem by constructing a third measure $\tilde{\mu}$ with $\tilde{\mu}\leq_{\mathrm{st}}\mu^*$ for all $(n_k)_{k\geq 1}$.  By showing that in fact $\tilde{\mu}=\mu^*$ we can then complete the proof by using features of both  $\mu^*$ and $\tilde{\mu}$ to compare $\mu^*|_{\La}$ with $\mu|_{\La}$.

Given $(s_x)_{x\in\GG}\in\Pi$ and a finite subgraph $\mathbf{H}\subseteq\GG$, Algorithm \ref{alg:PCF} gives us a function $\psi_{\mathbf{H}}$ taking $(s_x)_{x\in\GG}$ to a PCF process on $\mathbf{H}$, $(\eta^t_{\mathbf{H}})_{t\in[0,\infty]}$ say.  We now define $\tilde{\psi}_{\GG_n}$ to be the infimum of $\psi_{\mathbf{H}}$ taken over all finite $\mathbf{H}$ with $\GG_n\subseteq\mathbf{H}\subseteq\GG$.  Therefore
\begin{align}
\tilde{\psi}_{\GG_n}((s_x)_{x\in\GG})^t(e)&=\left\{\begin{array}{ll} 1 & \text{if }\eta^t_{\mathbf{H}}(e)=1 \text{ for all finite }\mathbf{H}\supseteq\GG_n \\ 0 & \text{otherwise}
\end{array} \right. , \\
\tilde{\psi}_{\GG_n}((s_x)_{x\in\GG})^t(v)&=\left\{\begin{array}{ll} w & \text{if }\eta^t_{\mathbf{H}}(v)=w \text{ for all finite }\mathbf{H}\supseteq\GG_n \\ f & \text{otherwise}
\end{array} \right. .
\end{align}
If we write $\tilde{\psi}_{\GG_n}((s_x)_{x\in\GG})=(\tilde{\eta}^t_{\GG_n})_{t\in[0,\infty]}$ then it is clear that $\tilde{\eta}^t_{\GG_n}\leq\eta^t_{\GG_n}$ for all $t$.  Moreover if $\GG_m\subseteq\GG_n$ then since $\{\mathbf{H}:\GG_m\subseteq\mathbf{H}\text{ finite}\}\supseteq\{\mathbf{H}:\GG_n\subseteq\mathbf{H}\text{ finite}\}$ we see that $\tilde{\eta}^t_{\GG_m}\leq\tilde{\eta}^t_{\GG_n}$ for all $t$.  Because of this monotonicity $(\tilde{\eta}^t_{\GG_n})_{t\in[0,\infty]}$ must converge to some process $(\tilde{\eta}^t)_{t\in[0,\infty]}$ as $n\longrightarrow\infty$.  Therefore if we define $\tilde{\mu}_{\GG_n,\alpha}$ to be the measure associated with $(\tilde{\eta}^t_{\GG_n})$ and $\tilde{\mu}$ to be the measure associated with $(\tilde{\eta}^t)$, then we see that the measures $\tilde{\mu}_{\GG_n,\alpha}$ converge to $\tilde{\mu}$.

It is easy to see from the construction that the limit $\tilde{\mu}$ does not depend on $\GG_1\subseteq\GG_2\subseteq\ldots$, and that $\tilde{\mu}$ is translation invariant.  This is because for each $e\in\mathcal{E}$ we have that
\begin{align*}
\tilde{\eta}^t(e)=0 \text{ }&\Leftrightarrow\text{ for all } n \text{ there exists a finite } \mathbf{H}\supseteq\GG_n \text{ such that } \eta_{\mathbf{H}}^t(e) = 0 \\
&\Leftrightarrow\text{ for all } n \text{ there exists a finite }  \mathbf{H}\supseteq\mathbf{H}_n \text{ such that } \eta_{\mathbf{H}}^t(e) = 0 ,
\end{align*}
where $\mathbf{H}_1\subseteq\mathbf{H}_2\subseteq\ldots$ is any other exhaustion of $\GG$ with $\mathbf{H}_n$ finite for all $n$.  Likewise $\tilde{\eta}^t(v)$ does not depend on the exhaustion for each $v\in\mathcal{V}$.  To see that translation invariance also holds let $\mathbf{H}_n=f(\GG_n)$ where $f$ is any isometry of $\GG$ and follow the argument of (\ref{eq:translation}).  Since $\tilde{\mu}$ is translation invariant the proof of Proposition \ref{prop:unique1} applies to $\tilde{\mu}$ and so we see that for all $t$ any warm infinite cluster of $\tilde{\eta}^t$ is unique almost surely.  Therefore by keeping $(s_x)_{x\in\GG}\in\Pi$ fixed we can use $(\tilde{\eta}^t)_{t\in[0\infty]}$ to define a coupled sequence of random freezing times $0<T_1<T_2<\ldots$ as follows.

\begin{definition}[Coupled freezing times]\label{def:coupledfreezing}$ _{}$

Let $(s_x)_{x\in\GG}\in\Pi$ be fixed and suppose $(\tilde{\eta}^t)_{t\in[0,\infty]}$ is as defined above.  For each time $t$ let $\ell(t)=\inf\{k:v_k \text{ is in a warm infinite cluster of }\tilde{\eta}^t\}$ with the convention that $\ell(t)$ is infinite if no warm infinite cluster exists at $t$.  Now suppose $X_{\infty,1},X_{\infty,2},\ldots$ are independent $\mathrm{Exp}(\alpha)$ random variables and define
\begin{align}
T_1 = \left\{\begin{array}{ll} X_{\infty,1} & \text{if } \ell(X_{\infty,1}) = \infty \\
\inf\{t:s_{v_{\ell(t)}}=t\} & \text{otherwise}\end{array}\right. .
\end{align}
Observe that if $\ell(X_{1,\infty})<\infty$ then $v_{\ell(X_{1,\infty})}$ must have joined the infinite cluster at time $t<X_{1,\infty}$.  Moreover because $v_{\ell(X_{1,\infty})}$ must have been warm at $t$ then $s_{v_{\ell(t)}}>t$ (and $s_{v_{\ell(t)}}-t\sim\mathrm{Exp}(\alpha)$).  Therefore $\{t:s_{v_{\ell(t)}}=t \}\neq\emptyset$ and so the memoryless property tells us that $T_1\sim\mathrm{Exp}(\alpha)$.

$T_2,T_3,\ldots$ are now defined inductively.  Suppose we have defined $T_1,\ldots,T_j$, then define
\begin{align}
T_{j+1} = \left\{\begin{array}{ll} T_j+X_{\infty,j+1} & \text{if } \ell(T_j+X_{\infty,j+1}) = \infty \\
\inf\{t:s_{v_{\ell(t)}}=t\} & \text{otherwise}\end{array}\right. ,
\end{align}
and note that we again have $T_{j+1}-T_j\sim\mathrm{Exp}(\alpha)$ from the memoryless property.
\end{definition}

Now consider the warm PCF process $(\zeta^t)_{t\in[0,\infty]}=\phi((s_x)_{x\in\GG})$, and let $\tau=\inf\{t:\zeta^t$ contains an infinite cluster$\}$.
If $T_1 = s_v$ for some $v$, then $T_1\geq\tau$ and $v$ is in an infinite cluster of $\tilde{\eta}$.  Lemma \ref{lem:dominate} tells us that $\tilde{\eta}\leq\zeta$ and so $v$ must also be in an infinite cluster of of $\zeta$ (which can not freeze).  Therefore conditioning on $(\eta^t)$ gives no further information about $s_v$.  Thus the law of $T_1$ given $(\eta^t)$ is
\begin{align}
\pr (T_1 \in . \,| (\eta^t)) &= \pr(T_1 \in . \,| (\eta^t);X_{\infty,1}<\tau) \times \pr(X_{\infty,1}<\tau|(\eta^t)) \nonumber \\
&\qquad + \pr(T_1 \in . \,| (\eta^t);X_{\infty,1}\geq \tau) \times \pr(X_{\infty,1}\geq \tau |(\eta^t)) \nonumber \\
&= \pr(X_{\infty,1} \in . \,| X_{\infty,1}<\tau) \times \pr(X_{\infty,1}<\tau) \label{eq:independent} \\
&\qquad + \pr(s_v \in . \,| s_v \geq \tau) \times \pr(X_{\infty,1}\geq \tau ) \nonumber\\
&= \pr(X_{\infty,1}\in .\,) \nonumber ,
\end{align}
and so $T_1$ is independent of $(\eta_t)$.

\begin{claim'}
For $t\in[0,T_1)$ we have $\tilde{\eta}^t=\zeta^t$.
\bpc
Consider the corresponding processes on $\GG_n$.  Lemma \ref{lem:monotonicity} and  Lemma \ref{lem:dominate} combine to tell us that if $C$ is an interior cluster of $\zeta_n^t$ then $C$ is also an interior cluster of $\zeta_{\mathbf{H}}^t$ for all finite $\mathbf{H}\supseteq\GG_n$, and so $C$ must also be a cluster of $\tilde{\eta}_n^t$.  Since $\GG_n\longrightarrow\GG$ we must therefore have that if $C$ is a finite cluster of $\zeta^t$ then $C$ is also a finite cluster of $\tilde{\eta}^t$.  Conversely if $C$ is an interior cluster of $\tilde{\eta}_n^t$ then $C$ is also an interior cluster of $\eta^t_n$ and $\zeta^t_n$.  Therefore any finite cluster of $\tilde{\eta}^t$ is also a finite cluster of $\zeta^t$.

Now suppose some edge $e\in\mathcal{E}$ has $\tilde{\eta}^t(e)\neq\zeta^t(e)$, then $e$ must be in an infinite cluster of $\zeta^t$ with $\zeta^t(e)=1$ and $\tilde{\eta}^t(e)=0$.  Since $\zeta^t(e)=1$ we must have $s_e<t$, therefore there is some time $s<s_e$ such that if $e=v_i v_j$ then either $v_i$ or $v_j$ has frozen in $\tilde{\eta}^s$ but not in $\zeta^s$.  Suppose it is $v_i$, then because we know that the finite clusters of $\tilde{\eta}^s$ and $\zeta^s$ are equal we must then have $v_i$ in an infinite cluster of $\tilde{\eta}^s$.  However, this is impossible since no infinite cluster freezes in either model before time $T_1$.  Thus $\tilde{\eta}^t=\zeta^t$ for all $t\in[0,T_1)$.
\epc
\end{claim'}

We now consider $(\zeta_{T_1}^t)_{t\in[0,\infty]}=\phi_{T_1}((s_x)_{x\in\GG})$.  We know that $\zeta_{T_1}^{T_1^-}=\tilde{\eta}^{T_1^-}$, and so because any infinite cluster of $\tilde{\eta}$ freezes at $T_1$ then we also have $\zeta_{T_1}^{T_1}=\tilde{\eta}^{T_1 }$.  By following the argument above and considering only the edges and vertices in $\GG\setminus\{\text{frozen infinite cluster}\}$ we can extend this to $\tilde{\eta}^t=\zeta_{T_1}^t$ for all $t\in[0,T_2)$.  Moreover since $T_2=T_1+X_{\infty,2}$ or $T_2=s_v$ for some $v$ in a warm infinite cluster of $\zeta_{T_1}$, we can repeat the argument of (\ref{eq:independent}) to see that $T_2$ is independent of $(\zeta_{T_1}^t)$.

By iterating we have that $\zeta_{T_1,\ldots,T_j}^{t}=\tilde{\eta}^{t}$ for all $t\in[0,T_{j+1})$.  Furthermore $T_j$ is independent of $(\zeta_{T_1,\ldots,T_k}^{t})$ for all $j>k$.  Therefore if we set $\mathcal{T}=\{0<T_1<T_2<\ldots\}$ with the $T_j$ defined by Definition \ref{def:coupledfreezing} then $\tilde{\eta}^t=\zeta_{\mathcal{T}}^{t}$ for all $t$.  Because $\mathcal{T}$ has independent $\mathrm{Exp}(\alpha)$ increments and is independent of $\tilde{\eta}^t$, then the law of $\zeta_{\mathcal{T}}^{t}$ must be $\mu$, and so $\tilde{\mu}=\mu^*$.

It now only remains to compare $\mu|_{\La}$ and $\mu^*|_{\La}$.  We do this by showing that for the coupled processes $(\eta_n^t)$ and $(\zeta_{\mathcal{T}}^{t})$ we have $\eta^t_n|_{\La}=\zeta_{\mathcal{T}}^{t}|_{\La}$ for all $t$ and all $n$ sufficiently large.  If this is true then it also follows that $\mu_n|_{\La}\longrightarrow\mu|_{\La}$ as required.

Suppose some edge or vertex $x\in\La$ is in a finite cluster of $\zeta_{\mathcal{T}}^{\infty}$, then $x$ is in some interior cluster of $\zeta_n^{\infty}$ for $n$ sufficiently large and so by Lemma \ref{lem:dominate} we know that $\eta_n^t(x)=\zeta_{\mathcal{T}}^{t}(x)$ for all $t$ and all $n$ sufficiently large.  Now suppose instead that $x\in\La$ is in an infinite cluster of $\tilde{\eta}^{\infty}=\zeta_{\mathcal{T}}^{\infty}$, and that this infinite cluster froze at time $T_k$.  Let $v=v_{\ell(T_k)}$ then at time $T_k -$ there must have been some warm path from $x$ to $v$ in $\tilde{\eta}^{T_k^-}$, and so for $n$ sufficiently large this path must also be present in $\tilde{\eta}^{T_k^-}_n$ and $\eta^{T_k^-}_n\geq\tilde{\eta}^{T_k^-}_n$.  Since $\eta_n^{T_k^-}\leq \zeta_{n,\mathcal{T}}^{T_k}$ we know that $\eta_n$ does not have $x$ joined to a vertex of higher priority than $v$, and thus $x$ will freeze in $\eta_n$ at time $T_k$.  Since the state of an edge $e$ is determined only by $s_e$ and the frrezing times of the adjacent vertices we must again have $\eta_n^t(x)=\zeta_{\mathcal{T}}^{t}(x)$ for all $t$ and all $n$ sufficiently large.  Therefore, since $\La$ is finite, when $n$ is sufficiently large we must have $\eta^t_n(x)=\zeta_{\mathcal{T}}^{t}(x)$ for all $t\in[0,\infty]$ and $x\in\La$.  This completes the proof.
\ep
\end{theorem}

\subsection{PCF on a infinite tree}

We have shown that there is an infinite volume PCF process on every amenable vertex transitive graph $\GG$.  The theorem below now tells us that an infinite volume PCF process also exist on any countably infinite tree $\mathbf{T}$.  Whether or not the PCF process can be defined on any graph remains open -- see Section \ref{sec:problems}.

\begin{theorem}
For every countable tree $\mathbf{T}=(\mathcal{V},\mathcal{E})$, and every fixed rate of freezing $\alpha>0$, there exists an infinite volume PCF process on $\mathbf{T}$.  This induces the PCF measure $\mu_{\mathbf{T},\alpha}$ on $\{0,1\}^{\mathcal{E}}$.
\bp
Given a (countably) infinite tree, $\mathbf{T}=(\mathcal{V},\mathcal{E})$, we can assign some vertex $v_0\in\mathcal{V}$ to be the root.  It is then possible to label the remaining vertices $\mathcal{V}\setminus\{v_0\}=\{v_1,v_2,\ldots\}$ in such a way that for each $v_j$ the unique path from $v_0$, $v_0 v_{i_1} \ldots v_{i_k}v_j$ say, has $i_1<\ldots<i_k<j$.  Suppose we have such a labelling and let $v_0\in\GG_1\subseteq\GG_2\subseteq\ldots$ be an exhaustion of $\mathbf{T}$ with each $\GG_n$ connected and finite.

By now assigning independent $\mathrm{Exp}(\alpha)$ random variables $X_v$ to each $v\in\mathcal{V}$ and independent $\mathrm{Exp}(1)$ random variables $X_e$ to each $e\in\mathcal{E}$, we can use the coupling of Definition \ref{def:coupling} to construct coupled processes on each of the $\GG_n$ -- $(\eta_n^t)_{t\in[0,\infty]}$ with measures $\mu_{n,\alpha}=(\mu_{n,\alpha}^t:t\in[0,\infty])$ say.  Now given $n$, the opening of any edge $e=v_i v_j\in\GG_n$ is determined only by $X_e$ and the states of $v_i$ and $v_j$.  Moreover, the state of a vertex $v\in\GG_n$ is determined by the $X_{v'}$ -- where $v'$ is the vertex of highest priority in the cluster containing $v$.  Note that our choice of enumeration means we must have $v'\in\{v_0,v_{i_1},\ldots,v_{i_k},v\}$ and so we need only know that states of the edges $v_0v_{i_1}, v_{i_1}v_{i_2},\ldots,v_{i_k}v$.  Since all of these edges and vertices must be in $\GG_n$ then if we consider PCF on $\GG_m\supseteq\GG_n$ for some $m>n$, then transitions of each of the edges and vertices of $\GG_n$ remain unchanged.  Therefore we see that $\eta_n^t|_{\GG_m}=\eta_m^t$ for all $t$, and so $\mu_{m,\alpha}|_{\GG_n}=\mu_{n,\alpha}$, implying that $\mu_{n,\alpha}$ converges as $n\longrightarrow\infty$.  The limit, $\mu_{\mathbf{T},\alpha}$, is the infinite volume rate $\alpha$ PCF measure on $\mathbf{T}$.
\ep
\end{theorem}

%% file: PCF_Zd.tex
\section{Proof of Proposition \ref{prop:no_infinite}}

On the cubic lattice $\ZZ^d=(\mathcal{V},\mathcal{E})$ we say a vertex $v\in\mathcal{V}$ is \emph{$\star$-open} in the PCF model if at least one of its incident edges is open,  otherwise we say that $v$ is \emph{$\star$-closed}.  Observe that if a PCF configuration $\eta$ contains an infinite cluster of open edges then it must also contain an infinite $\star$-open cluster of vertices.

Now consider the site percolation model on $\ZZ^d$ (where vertices are open with probability $p$ and closed with probability $1-p$ independently of each other).  It is well know that there exists a critical value $0<p_c^{\mathrm{site}}<1$ such that if $p<p_c^{\mathrm{site}}$ then all open clusters are almost surely finite -- see e.g. \cite{grimmett1} for details.  Therefore to prove Proposition \ref{prop:no_infinite} it suffices to show that when $\alpha>0$ is sufficiently large the PCF measure of $\star$-open vertices is stochastically dominated by some site percolation measure $\pr_p$ with $p<p_c^{\mathrm{site}}$.  This will be a consequence of the lemma below.

\begin{lemma}
Set $d\geq 2$, and consider rate $\alpha$ PCF on $\ZZ^d$.  Let $v\in\mathcal{V}$, and $A$ be any event which is measurable with respect to the $\star$-configuration of $\mathcal{V}\setminus\{v\}$, then
\begin{align*}
\mu_{\ZZ^d,\alpha}(v \text{ is } \star\text{-open}\,|\,A) \leq f_d(\alpha) ,
\end{align*}
where $f(\alpha)\longrightarrow 0$ as $\alpha\longrightarrow\infty$.
\bp
We use a coupling argument.  To each vertex $v\in\mathcal{V}$ we associate an independent $\mathrm{Exp}(\alpha)$ random variable $X_v$ and $2d$ independent $\Gamma(\frac{1}{2},1)$ random variables ($Y_{v,1}\ldots Y_{v,2d}$).  Now at each edge $e=v_1 v_2\in\mathcal{E}$ choose unique $Y_{v_1,i}$ and $Y_{v_2,j}$ from $v_1$ and $v_2$ and set $X_e = Y_{v_1,i}+Y_{v_2,j}$.  Note that since $X_e$ is the sum of two independent $\Gamma(\frac{1}{2},1)$ random variables we have $X_e\sim\mathrm{Exp}(1)$.  We can now use these random variables to drive PCF on $\ZZ^d$ (as in Algorithm \ref{alg:PCF}).  Note that the freezing of any infinite clusters will again be controlled by independent exponential clocks; but since a $\star$-open vertex is never in an infinite cluster this is not a problem.

Now a vertex $v$ is $\star$-open in the PCF model only if an adjoining edge opens before $v$ freezes.  For this we require
\[ X_v>\mathrm{min}\{X_e:e\text{ is adjacent to }v\}>\mathrm{min}\{Y_{v,1}\ldots Y_{v,2d}\}.
\]
Therefore we are done by setting
\[ f_d(\alpha) = \pr (X_v>\mathrm{min}\{Y_{v,1}\ldots Y_{v,2d}\}) \longrightarrow 0 \quad \text{as} \quad \alpha \longrightarrow \infty . \qedhere
\]
\ep
\end{lemma}

Now if we set $\alpha^* = \mathrm{sup}\{\alpha  :  f_d(\alpha) \geq p_c(d) \}$ then for $\alpha>\alpha^*$ the measure of $\star$-open vertices is stochastically dominated by $\pr_p$ with $p<p_c$ and so all the clusters in the PCF configuration $\eta$ on $\ZZ^d$ are finite almost surely.

Observe that this proof only relies on $\ZZ^d$ having $p_c^{\mathrm{site}}>0$, and $\ZZ^d$ having maximum degree $2d<\infty$.  Therefore this proof allows Proposition \ref{prop:no_infinite} to be extended to any amenable vertex transitive graph $\LL^d$ where the degree of every vertex is bounded by some $\Delta<\infty$.

\begin{remark}
For $d=2$ we have $p_c^{\mathrm{site}}\approx 0.59$ and so a back of the envelope calculation gives $\alpha^*\approx 20$.  This upper bound for $\alpha_c$ far exceeds the estimate of $\alpha_c\approx 0.55$ found in Section \ref{sec:sim}.
\end{remark}

%% file: PCF_tree.tex
\section{The results on a tree}\label{sec:tree}

This section begins with a restatement of Theorem \ref{thm:cluster_size} in greater generality.

\begin{theorem}\label{thm:treetail}
Let $\alpha>0$ be fixed and consider a rate $\alpha$ PCF process on the rooted $d$-ary tree $\mathbf{T}_d$ -- with measure $\mu_{\mathbf{T}_d,\alpha}$.  Let $A_k$ be the event that the root cluster has size $k$ ($k$ vertices) in the final distribution.  Then
\vspace{-0.5\baselineskip}\begin{itemize}
\item if $\alpha<d-1$ we have $\mu_{\mathbf{T}_d,\alpha}(A_k)\sim C_{d,\alpha}\, k^{-2}$, for some fixed constant $C_{d,\alpha}$.
\item If $\alpha=d-1$ we have $\mu_{\mathbf{T}_d,\alpha}(A_k)\sim C_{d}\, k^{-\left(2-\frac{1}{2d}\right)}$, for some fixed constant $C_d$.
\item If $\alpha>d-1$ then $\mu_{\mathbf{T}_d,\alpha}(A_k)$ decays exponentially in $k$.
\end{itemize}
\end{theorem}

Before embarking upon the proof of Theorem \ref{thm:treetail} we remark that this result is also valid for unrooted $d$-ary trees.  Moreover, with only minor modifications to the proof one can show that the result can be extended to and tree $\mathbf{T}$ where the limit
\begin{align*}
d = \lim_{n\rightarrow\infty}\sqrt[n]{|\partial\mathbf{\Lambda_n}|}
\end{align*}
exists.  Here $|\partial\Lambda_n|$ gives the number of vertices at distance $n$ from the root.  In this case the $d$ plays the same role as in Theorem \ref{thm:treetail}.

\subsection{Time rescaling}

As mentioned in the introduction our understanding of PCF on the tree comes through viewing it as time rescaled percolation.  The following proposition makes sense of this.

\begin{proposition}\label{prop:timerescale}
Suppose we have a rate $\alpha$ PCF process on a tree $\mathbf{T}$.  Let $W_v(t)$ be the event that $\eta_t(v)=w$ -- i.e. that the vertex $v$ is still warm at time $t$.  Suppose $C$ is any finite connected component containing $v$, then if we condition on $W_v(t)$ the probabilities of $\{C_v(t) \subseteq C\}$ and $\{C_v(t) = C\}$ are given by
\begin{align}
\mu^t_{\mathbf{T},\alpha} (\{C_v(t) \subseteq C\} | W_v(t) ) &= p^{| C|} \label{eq:PCFtreeprob1}, \\
\text{and} \quad \mu^t_{\mathbf{T},\alpha} (\{C_v(t) = C\} | W_v(t) ) &= p^{|C|} (1-p)^{|\partial C|} \label{eq:PCFtreeprob2}.
\end{align}
Where
\begin{align}
p=\frac{1-\ee^{-(1+\alpha)t}}{1+\alpha} \label{eq:PCFtreeprob3}.
\end{align}
\bp
Label the edges of $C$, $e_1,e_2,\ldots,e_{|C|}$ in such a way that the labels of any path away from $v$ are increasing. Then
\begin{align}
\mu^t_{\mathbf{T},\alpha} \left( \{C_v(t) \subseteq C \}\,|\, W_v(t) \right) = \,& \mu^t_{\mathbf{T},\alpha} \left( e_1 \text{ open at } t \,|\, v \text{ warm} \right) \times
\mu^t_{\mathbf{T},\alpha} \left( e_2 \text{ open at } t \,|\, v, e_1 \text{ warm} \right) \nonumber \\ & \times \cdots \times
\mu^t_{\mathbf{T},\alpha} \left( e_{| C|} \text{ open at } t \,|\, v, e_1, \ldots, e_{| C| -1} \text{ warm} \right) \nonumber , 
\end{align}
and
\begin{align}
\mu^t_{\mathbf{T},\alpha} \left( \{C_v(t) = C \}\,|\, W_v(t) \right) = \,& \mu^t_{\mathbf{T},\alpha} \left( C_v(t) \subseteq C \right) \times \mu^t_{\mathbf{T},\alpha} \left( \text{boundary edges closed} \,|\, C_v(t) \text{ warm}\right) \nonumber .
\end{align}

Our choice of labelling means that if $v$ is warm at $t$ and $e_1,\ldots,e_{j-1}$ are open then the end of $e_j$ closest to $v$ must also be warm at time $t$.  Thus
\begin{align}
\mu^t_{\mathbf{T},\alpha} \left( e_{j} \text{ open at } t \,|\, v, e_1, \ldots, e_{j-1} \text{ warm} \right) &= \mu^t_{\mathbf{T},\alpha} \left( e_j \text{ opens before } t \text { and before its other end freezes} \right) \nonumber \\
&= \int_{0}^{t} \left(1-\ee^{-\alpha t}\right) \ee^{-t}\, \dd t \nonumber \\
&= \frac{1-\ee^{-(1+\alpha)t}}{1+\alpha} \label{eq:rescaleproof2}.
\end{align}
Since we are working on a tree each of the boundary edges is independent and so we get
\begin{align}
\mu^t_{\mathbf{T},\alpha} \left( \text{boundary edges closed} \,|\, C_v(t)\text{ warm} \right) &= \mu^t_{\mathbf{T},\alpha} \left( \text{boundary edge closed} \,|\, C_v(t)\text{ warm} \right)^{|\partial C|} \nonumber ,
\end{align}
where
\begin{align}
\mu^t_{\mathbf{T},\alpha} \left( \text{boundary edge closed} \,|\, C_v(t)\text{ warm} \right) &= 1 - \frac{1-\ee^{-(1+\alpha)t}}{1+\alpha} = \frac{\alpha + \ee^{-(1+\alpha)t}}{1+\alpha} .
\end{align}
The result then follows.
\ep
\end{proposition}

To understand the sense in which Proposition \ref{prop:timerescale} gives a time rescaling, observe that if $\alpha=0$ then we have a bond percolation process where edges open independently at rate 1.  Let $\pr_{\tau}$ denote the measure of this process. In this case the probability that a given edge is open at time $\tau$ is $\pr_{\tau}(e\text{ open at }\tau)=1-\ee^{-\tau}$.  Comparing this with (\ref{eq:PCFtreeprob3}) we see that
\begin{align*}
\mu^t_{\mathbf{T},\alpha} (\{C_v(t) \subseteq C\} | W_v(t) ) &= \pr_{\tau} (\{C_v(\tau) \subseteq C\} ), \\
\text{and} \quad \mu^t_{\mathbf{T},\alpha} (\{C_v(t) = C\} | W_v(t) ) &= \pr_{\tau}(\{C_v(\tau) = C\}) ,
\end{align*}
where
\begin{align}
\tau = - \log \left( \frac{\alpha+\ee^{-(1+\alpha)t}}{1+\alpha} \right) \label{eq:time rescale}.
\end{align}

Note that the distribution of infinite components is entirely determined by the events $C_v(t) \subseteq C$ for finite $C$.  Therefore the time rescaling of Proposition \ref{prop:timerescale} applies equally to finite and infinite clusters.  Because there is always a positive probability that a given cluster $C_v(t)$ is warm, we deduce from the proposition that running PCF on $\mathbf{T}$ to time $t$ yields infinite clusters if and only if the bond percolation on $\mathbf{T}$ yields infinite clusters at time $\tau$ (where $\tau$ is given by (\ref{eq:time rescale})).

%
%

\begin{remark}\label{rem:Knlimitingcase}
In \cite{bn&k1} and \cite{bn&k2}, Ben-Naim and Krapivsky work on the complete graph $\mathbf{K}_n$, opening edges at a rate $\frac{1}{n}$ and allowing vertices to freeze at rate $\alpha$.  Using a Smoluchowski equation they obtain a result which is similar to Proposition \ref{prop:timerescale} showing that in this case -- conditional on a cluster being warm -- the cluster's size at time $t$ is equal in distribution to that of a percolation cluster at time 
\[ \tau = \frac{1-\ee^{-\alpha t}}{\alpha} .
\]
In fact it is possible for us to deduce this result from Proposition \ref{prop:timerescale}.  We start with the fact that for large $n$ small clusters are cycle free with high probability, and therefore for $0\leq t< 1$ they grow in the same way as cluster on an $n$-ary tree.  See chapter 2 of \cite{durrett} for details.  Now because we are opening edges on $\mathbf{K}_n$ at rate $\dfrac{1}{n}$, we need to rescale time to $\dfrac{t}{n}$ and the rate of freezing to $n\alpha$.  Putting this in equation (\ref{eq:PCFtreeprob3}) we get
\begin{align*}
\tau = np = n \, \frac{1-\ee^{-(1+n\alpha)\frac{t}{n}}}{1+n\alpha} 
\longrightarrow \frac{1-\ee^{-\alpha t}}{\alpha} \quad \text{as} \quad n\longrightarrow\infty .
\end{align*}
Therefore we see that for finite clusters, mean field PCF can be thought of as the limiting case of PCF on a $n$-ary tree as $n\longrightarrow \infty$.
\end{remark}

\subsection{The critical phenomenon}

Consider an infinite tree $\mathbf{T}$ with a distinguished vertex $v$, suppose also that $\mathbf{T}$ satisfies the  \emph{exponential volume growth property}
\begin{align}
d = \liminf_{n \rightarrow \infty} \sqrt[n]{\left| \Gamma^{n}(v) \right|}  > 1 \label{eq:expgrowth} ,
\end{align}
where $\left| \Gamma^{n}(v) \right|$ is the number of vertices at distance $n$ from $v$.  Then if we run percolation (without freezing) on $\mathbf{T}$ until time $\tau$, there is a strictly positive probability that $C_v(\tau)$ is infinite ($\pr_{\tau}(\{|C_v(\tau)|=\infty\})>0$) if and only if \begin{align*}
\pr(e \text{ open at }\tau) = 1-\ee^{-\tau} > \frac{1}{d} .
\end{align*}
This follows from considering $C_v(\tau)$ as having grown from a branching process.  For details see e.g. Section 2.1 of \cite{durrett}.

\begin{proposition}
Consider rate $\alpha$ PCF running on an infinite tree $\mathbf{T}$ with distinguished vertex $v$.  Suppose $\mathbf{T}$ has growth rate $d = \liminf_{n \rightarrow \infty} \sqrt[n]{\left| \Gamma^{n}(v) \right|}  > 1$.  Then the critical rate of freezing for the existence of infinite clusters is $\alpha_c = d-1$.  That is
\begin{align}
\mu_{\mathbf{T},\alpha}(\{|C_v(\infty)|=\infty\}) \left\{\begin{array}{ll} =0 & \text{for } \alpha \geq \alpha_c \\ >0 & \text{for } \alpha < \alpha_c \end{array} \right. .
\end{align}
Moreover, in the case $\alpha < \alpha_c$, the critical time for the emergence of infinite clusters is
\begin{align}
t_c=\dfrac{1}{1+\alpha}\log\left(\dfrac{d}{1+\alpha}\right) .
\end{align}
\bp
From Proposition \ref{prop:timerescale} we know that if we condition on $W_v(t)$, then the distribution of the structure of $C_v(t)$ is equal to that of percolation cluster where edges are open with probability $p=\frac{1}{1+\alpha} \left( 1-\ee^{-(1+\alpha)t} \right)$.  Let $\pr_p(\{|C_v|=\infty\})$ be the probability that the cluster containing $v$ is infinite in the percolation model.  Then given $\alpha < \alpha_c = d-1$ for all $t>t_c$ we have $p>\frac{1}{d}$, implying $\pr_p(\{|C_v|=\infty\})>0$.  Therefore
\begin{align}
0 < \mu^t_{\mathbf{T},\alpha}(W_v(t)) \times \pr_p(\{|C_v|=\infty\}) \leq \mu^t_{\mathbf{T},\alpha}(\{|C_v(t)|=\infty \}) \leq \mu_{\mathbf{T},\alpha}(\{|C_v|=\infty\} ),
\end{align}
since $\mu^t_{\mathbf{T},\alpha}(W_v(t)) = \ee^{-\alpha t} > 0$ for all $t$.

However, if $t<t_c$, then $p<\frac{1}{d}$ and thus $\mu^t_{\mathbf{T},\alpha}(\{|C_v(t)|=\infty\}) \leq \pr_p(\{|C_v|=\infty\}) =0$.  Hence the time $t_c$ is critical.

Conversely suppose $\alpha \geq \alpha_c$.  If $C_v$ is infinite in the final distribution, then there must be some finite time $t$ at which $C_v(t)$ is infinite.  But
\begin{align}
\mu^t_{\mathbf{T},\alpha}(\{|C_v(t)|=\infty\}) \leq \pr_p(\{|C_v|=\infty\}) =0 \quad \forall \, t,
\end{align}
and so $\mu_{\mathbf{T},\alpha}(\{|C_v|=\infty\} )=0$.
\ep
\end{proposition}

From now on we shall restrict our focus to $d$-ary trees $\mathbf{T}_d$.  Observe that a $d$-ary tree has $\left| \Gamma^n(v) \right|=(d+1)d^{n-1}$, and a rooted $d$-ary tree has $\left| \Gamma^n(v) \right|=(d+1)d^{n-1}$.  Thus they both have exponential growth at rate $d$.

\begin{corollary}\label{thm:treecriticalvalue}
Suppose we run rate $\alpha$ PCF on the (possibly rooted) $d$-ary tree $\mathbf{T}_d$.  Then
\vspace{-0.5\baselineskip}\begin{itemize}
\item if $\alpha <  d-1$, then at time $t>t_c=\dfrac{1}{1+\alpha}\log\left(\dfrac{d}{1+\alpha}\right)$, $\mathbf{T}_d$ will contain infinite components almost surely.
\item if $\alpha \geq d-1$, then $\mathbf{T}_d$ will have no infinite component with probability $1$.
\end{itemize}
\end{corollary}

\subsection{The final distribution -- proof of Theorem \ref{thm:treetail}}

Heuristically, the reason that we have clusters of algebraic size in the super-critical regime is because our system must pass though a critical time
\begin{align}
t_c=-\dfrac{\log \left(1-\frac{1+\alpha}{d}\right)}{1+\alpha} ,
\end{align}
at which large clusters exist, but have not yet merged to become infinite.  Some of these large clusters can then freeze at or near this critical time leaving us with large finite clusters in the final distribution.  However, the process is only sufficiently close to this critical time for a short period, meaning that the exponent of the cluster sizes is larger than the value of $\frac{3}{2}$ we see for clusters in critical bond percolation.

In the case of critical PCF the system spends a much more significant amount of time approaching criticality, but the clusters never become infinite.  Therefore there are more large clusters to freeze, thus explaining why we have a lower exponent than in the super-critical case.  To see why the exponent should depend on the dimension we look at the time rescaling of equation (\ref{eq:time rescale}), and the rate at which it reaches its maximum value
\[ \tau_{\text{max}} - \tau = \frac{1}{1+\alpha} - \frac{1-\ee^{-(1+\alpha)t}}{1+\alpha} = \frac{\ee^{-(1+\alpha)t}}{1+\alpha} .
\]
Noting that $d=1+\alpha_c$ we see that the degree has a dramatic difference on the rate at which the system approach criticality, and therefore it is not surprising that the algebraic exponent depends on the dimension.

We now set about proving Theorem \ref{thm:treetail} formally by calculating the probabilities $\mu_{\mathbf{T}_d,\alpha}(A_k)$ explicitly.  We can do this by integrating the probability $\mu_{\mathbf{T}_d,\alpha}^{t}(\{C_v(t) = C\}\,|\,W_v(t))$ with respect to the time at which $v$ freezes -- $\alpha\ee^{-\alpha t}\dd t$.  Thus we get
\begin{align*}
\mu_{\mathbf{T},\alpha}(\{C_v(\infty) = C\}) &= \int_{0}^{\infty} \mu_{\mathbf{T},\alpha}^{t}(\{C_v(t)=C\}\,|\,W_v(t)) \, \alpha\ee^{-\alpha t}\dd t \nonumber \\
&= \int_{0}^{\infty} p^{|C|} (1-p)^{|\partial C|} \,\alpha\ee^{-\alpha t}\dd t ,
\end{align*}
where $p=\frac{1}{1+\alpha}\left(1+\ee^{-(1+\alpha)t}\right)$.  By rewriting this as an integral in terms of $p$ we get
\begin{align}
\mu_{\mathbf{T},\alpha}(\{C_v(\infty) = C\}) 
&= \alpha \int_{0}^{\frac{1}{1+\alpha}} p^{|C|} (1-p)^{| \partial C|} (1-(1+\alpha)p)^{-\frac{1}{1+\alpha}} \,\dd p \label{eq:treefinal2} \\
&= \frac{\alpha}{1+\alpha} \frac{1}{(1+\alpha)^{| C|}} \frac{\Gamma\left(|C|+1\right) \Gamma\left(\frac{\alpha}{1+\alpha}\right)}{\Gamma\left(|C| + \frac{1+2\alpha}{1+\alpha}\right) } \, {}_2F_1\left[\begin{matrix} | C| +1, -| \partial C| \\ | C| + \frac{1+2\alpha}{1+\alpha} \end{matrix}; \frac{1}{1+\alpha} \right] \label{eq:treefinal3} .
\end{align}
Here the last line follows from substituting $\tilde{p}=(1+\alpha)p$ and using Euler's hypergeometric transform.  Observe also that the probability $\mu_{\mathbf{T},\alpha}(\{C_v(\infty) \subseteq C\})$ can be calculated using the same formulae by setting $| \partial C| = 0$.

Whilst it is nice to have a closed form for $\mu_{\mathbf{T},\alpha}(\{C_v = C\})$ -- equation (\ref{eq:treefinal3}) -- the hypergeometric function can often prove impenetrable, and so for practical purposes we shall use (\ref{eq:treefinal2}).

To give explicit results on the rooted $d$-ary tree $\mathbf{T}_d$ we shall need the following two lemmas.  We have chosen to focus on the rooted tree in order to simplify calculations (especially Lemma \ref{lem:noktrees}).  However, the corresponding results for the unrooted case prove to be very similar, and the main result of this section, Theorem \ref{thm:treetail}, holds for both cases.

\begin{lemma}
Consider the rooted $d$-ary tree.  Suppose $C\subseteq\mathbf{T}_d$ is a connected component of size $k$ -- i.e. with $k$ vertices.  Then $|C| = k-1$ and $| \partial C| = (d-1)k+1$, and thus
\begin{align}
\mu_{\mathbf{T},\alpha}(\{C_v (\infty) = C\}) 
&= \alpha \int_{0}^{\frac{1}{1+\alpha}} p^{k-1} (1-p)^{(d-1)k+1} (1-(1+\alpha)p)^{-\frac{1}{1+\alpha}} \,\dd p \label{eq:dtreefinal} .
\end{align}
\end{lemma}

\begin{lemma}\label{lem:noktrees}
The number of trees on $k$ vertices which are rooted sub-graphs of the $d$-ary tree $\mathbf{T}^d$ is given by
\begin{align*}
N_k^d = \frac{1}{k}\binom{d\,k}{k-1} .
\end{align*}
\end{lemma}


The proofs of these lemmas are mainly combinatorial.  For details see e.g. \cite{klarner}.

\bp
We can now give a formal proof of Theorem \ref{thm:treetail}.  Throughout this proof we shall use $C_{d,\alpha}$, $C_d$ and $C_{\alpha}$ to represent constants whose values can change from line to line.  It is possible to calculate their values explicitly, but not enlightening to do so.
Using Stirling's formula one can show
\begin{align}
N^d_k=\frac{1}{k}\binom{d\,k}{k-1}\sim C_d \left(\frac{d^d}{(d-1)^{d-1}}\right)^k\frac{1}{k^{\frac{3}{2}}}\label{eq:dcat}. 
\end{align}
Therefore by combining equations (\ref{eq:dtreefinal}) and (\ref{eq:dcat}) we get
\begin{align}
\mu_{\mathbf{T}_d,\alpha}(A_k) &= \frac{1}{k}\binom{d\,k}{k-1}  \alpha \int_{0}^{\frac{1}{1+\alpha}} p^{k-1}  (1-p)^{(d-1)k+1}  \left( 1-(1+\alpha) p \right)^{- \frac{1}{1+\alpha}} \, \dd p \nonumber \\
&\sim C_{d,\alpha} \frac{1}{k^{\frac{3}{2}}} \int_{0}^{\frac{1}{1+\alpha}} f^{k-1}(x) g_{\alpha}(x) \, \dd x \label{eq:tailest}
\end{align}
where $f(x) = \frac{d^d}{(d-1)^{d-1}}x(1-x)^{d-1}$, $g_{\alpha}(x) = \frac{(1-x)^d}{\left(1-(1+\alpha) x \right)^{\frac{1}{1+\alpha}}}$.

Observe now that $f(x)$ is positive on $[0,1]$, and attains a unique maximal value of $1$ at $x=\frac{1}{d}$.
We also note that $g_{\alpha}$ is integrable on $[0,1]$ (for $\alpha>0$) and is continuous everywhere except $x=\frac{1}{1+\alpha}$.  We can now consider the different cases:

\paragraph*{\underline{$\alpha>d-1:$}}
Bounding $f(x)$ by its supremum we have
\begin{align}
\mu_{\mathbf{T}_d,\alpha}(A_k) &\lesssim C_{d,\alpha} \frac{1}{k^{\frac{3}{2}}}\int_{0}^{\frac{1}{1+\alpha}} \sup_{x\in\left[0,\frac{1}{1+\alpha}\right]}f^{k-1}(x) g_{\alpha}(x)\,\dd x \nonumber \\
&=C_{d,\alpha}\frac{1}{k^{\frac{3}{2}}}\left(\sup_{x\in\left[0,\frac{1}{1+\alpha}\right]} f(x) \right)^{k-1} \int_{0}^{\frac{1}{1+\alpha}}g_{\alpha}(x)\,\dd x ,
\end{align}
which decays exponentially in $k$ since $\sup_{x\in\left[0,\frac{1}{1+\alpha}\right]} f(x)<1$.

\paragraph*{\underline{$\alpha<d-1:$}}
It is relatively easy to verify that 
\begin{align*}
\mu_{\mathbf{T}_d,\alpha}(A_k)\sim C_{d,\alpha}\frac{1}{k^{\frac{3}{2}}} \left[ \int_0^1 f^{k-1}(x)\,\dd x + \mathcal{O}\left(\frac{1}{k}\right)\right].
\end{align*}
This integral can now be evaluated explicitly in terms of the Beta function as
\begin{align}
\int_0^1 f^{k-1}(x)\,\dd x = \left( \frac{d^d}{(d-1)^{d-1}} \right)^{k-1} \mathrm{B}(k,(d-1)k-d+2)\sim C_d \frac{1}{k^{\frac{1}{2}}} .
\end{align}
From which the result follows.

\paragraph*{\underline{$\alpha=d-1:$}}
In this case the singularity of $g_{\alpha}$ coincides with the maximum of $f$, and so its contribution becomes significant.  By expanding $f$ about $x=\frac{1}{d}$ we get
\begin{align}
f(x)=1-\frac{d^3}{2(d-1)}\left(x-\tfrac{1}{d}\right)^2 + \mathcal{O}\left(\left(x-\tfrac{1}{d}\right)^3\right) = \tilde{f}(x)+\mathcal{O}\left(\left(x-\tfrac{1}{d}\right)^3\right).
\end{align}
Therefore 
\begin{align}
\mu_{\mathbf{T}_d,\alpha}(A_k) &\sim C_d\left[ \frac{1}{k^{\frac{3}{2}}} \int_0^{\frac{1}{d}} \tilde{f}^{k-1}(x)\frac{(1-x)^d}{(1-dx)^{\frac{1}{d}}} \, \dd x +\mathcal{O}\left(\frac{1}{k}\right)\right] \nonumber \\
&= C_d\left[ \frac{1}{k^{\frac{3}{2}}} \left(1+\mathcal{O}\left(\frac{1}{\sqrt{k}}\right)\right)\int_0^{1} (1-y)^{k-1} y^{-\frac{1}{2d}} \, \frac{\dd y}{\sqrt{y}} +\mathcal{O}\left(\frac{1}{k}\right)\right] .
\end{align}
Where the second line follows from setting $y=1-\tilde{f}(x)$ and approximating $(1-x)^d$ by $\left(1-\frac{1}{d}\right)^d$ -- its value at $x=\frac{1}{d}$.  The integral can then be calculated explicitly as
\[\int_0^{1} (1-y)^{k-1} y^{-\frac{1}{2d}} \, \frac{\dd y}{\sqrt{y}} = \mathrm{B}\left(k,\frac{1}{2}-\frac{1}{2d}\right)\sim C_d \frac{1}{k^{\frac{1}{2}-\frac{1}{2d}}},
\]
and so the result follows.
\ep




%% file: PCF_simulations.tex
\section{Simulations of PCF on $\ZZ^2$}\label{sec:sim}

Using Algorithm \ref{alg:PCF} we can write a fairly efficient program to simulate PCF on a finite graph.  In this section we use Monte Carlo simulations of PCF on an $n$ by $m$ square grid in order to give us insight into PCF on $\ZZ^2$.  We begin by using simulations to estimate the critical value $\alpha_c$.

\subsection{Estimating $\alpha_c$}\label{sec:estimatealpha_c}

It is known (e.g. see \cite{grimmett1}) that for bond percolation on an $n$ by $n+1$ square grid, the probability that there exists a left--right crossing, $C_{LR}$, satisfies
\begin{align*}
\pr_p(C_{LR}) \left\{ \begin{array}{ll}
= \dfrac{1}{2} & \text{for }p=p_c=\frac{1}{2} \\
\displaystyle{\xrightarrow[n\rightarrow\infty]{} 0 }& \text{for }p<p_c=\frac{1}{2} \\
\displaystyle{\xrightarrow[n\rightarrow\infty]{} 1 }& \text{for }p>p_c=\frac{1}{2}
\end{array} \right. .
\end{align*}

\begin{figure}[ht]
\centering
\includegraphics[width=0.9\textwidth]{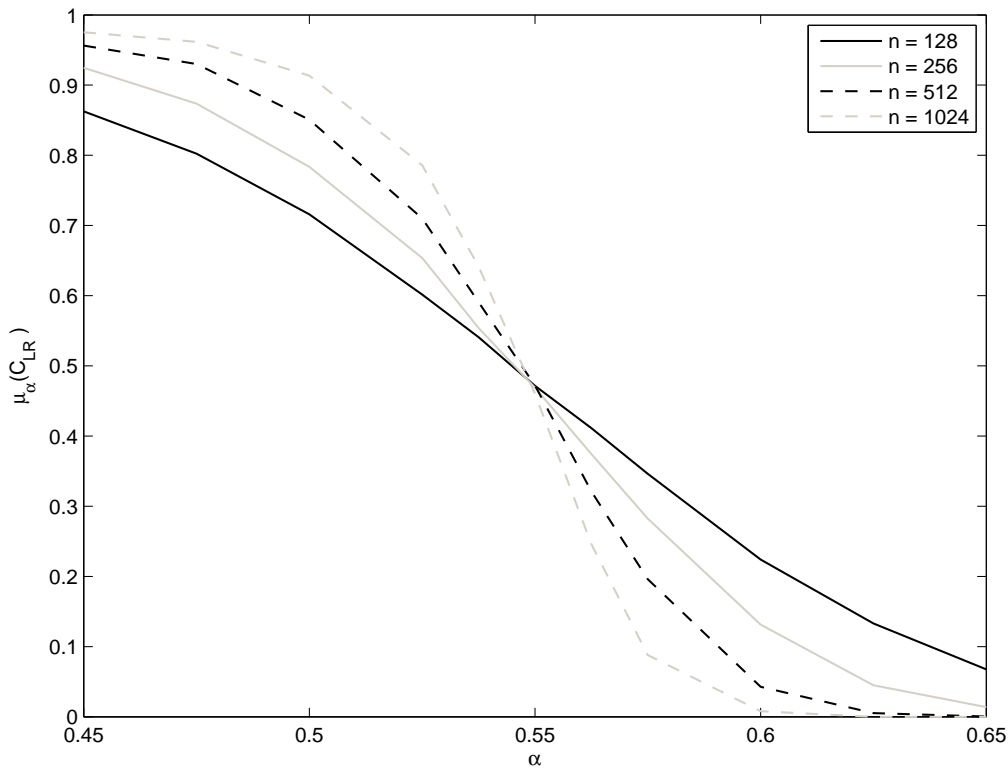}
\caption{Monte Carlo simulations for the probability of a left--right crossing of a $n$ by $n+1$ square grid when $n=128$, $256$, $512$ and $1024$.  For each $n$ at least 2500 simulations have been used for each data point.}\label{fig:left--right}
\end{figure}

Therefore, the existence of a left--right crossing of a $n$ by $n+1$ grid would seem to be a reasonable indicator for the existence of an infinite component for PCF on $\ZZ^2$.  Figure \ref{fig:left--right} shows a Monte Carlo simulation of how the probability of a left--right crossing varies with $\alpha$ for various values of $n$.  We focus on values of $\alpha$ between $0.45$ and $0.65$ since this is where the transition between sub-critical and super-critical PCF occurs.

Observe that for each $n$, we get $\mu_{\alpha}(C_{LR})=\frac{1}{2}$ when $\alpha \approx 0.55$, and therefore $\alpha_c \approx 0.55$ would seem to be a reasonable estimate for the critical rate of freezing on $\ZZ^2$.  Moreover, we see that as $n$ increases the transition becomes sharper -- i.e. changes in $\alpha$ have a greater effect on the existence of a left--right crossing, as we would expect.

\subsection{The cluster size distribution}\label{sec:tailsim}

It is known (for a reference see e.g. \cite{grimmett1}) that in sub-critical bond percolation on $\ZZ^d$ the distribution of the size of the clusters decays exponentially, and in super-critical bond percolation on $\ZZ^d$ the size of the finite clusters decays sub-exponentially.  This means to say that for $0< p<p_c$ or $p_c<p<1$ there exits a $\lambda(p)>0$ such that 
\begin{align*}
\pr_p(\{|C_0|=k\})\prec \left\{ \begin{array}{ll} \ee^{-\lambda(p)k} & \text{if } p<p_c \\ \ee^{-\lambda(p)k^{\frac{d-1}{d}}} & \text{if } p > p_c \end{array} \right. .
\end{align*}
Thus it is only at criticality, $p=p_c$, that large finite clusters are likely to exist.  In this case we have 
\begin{align*}
\pr_{p_c}(\{|C_0|=k\}) \asymp  k^{-\gamma(d)} ,
\end{align*}
for some constant $\gamma(d)$ which depends only on the dimension.
It is thought that in $2$ dimensions the exponent is $\gamma(2)=\frac{96}{91}\approx 1.055$, and that $\gamma(d)$ attains a mean field value of $\frac{3}{2}$ for dimension $d>6$ (see \cite{grimmett3}).

A similar result is true for bond percolation on a $d$-ary tree where
\begin{align*}
\pr_p(\{|C_0|=k\}) \left\{ \begin{array}{ll} \prec  \ee^{-\lambda(p)k} & \text{for } p \neq p_c \\ \asymp  k^{-\frac{3}{2}} & \text{for } p=p_c \end{array} \right. .
\end{align*}

However, as we discovered in Section \ref{sec:tree}, running PCF on a $d$-ary tree results in somewhat different behaviour, and we have a power law distribution for the cluster size in both the critical and super-critical regimes.  Indeed Theorem \ref{thm:treetail} tells us that
\begin{align*}
\mu_{\mathbf{T}_d,\alpha}(\{|C_0|=k\}) \left\{ \begin{array}{ll} \prec  \ee^{-\lambda(\alpha)k} & \text{for } \alpha > \alpha_c  \\ \asymp  k^{-\left(2-\frac{1}{2d}\right)} & \text{for } \alpha=\alpha_c \\ \asymp  k^{-2} & \text{for } \alpha < \alpha_c  \end{array} \right. ,
\end{align*}
also demonstrating that the exponents in the critical and super-critical regimes differ.

This happens because in a super-critical regime the process must pass through some critical time $t_c$ (depending on $\alpha<\alpha_c$) at which an infinite cluster will first appear.  Since it is possible for large clusters to freeze at this time rather than merge into an infinite cluster we are left with large clusters in our final distribution.  We account for the larger exponent in the super-critical case from the fact that the process only briefly passes through criticality, whereas for $\alpha=\alpha_c$ the length of time the process spends near criticality is much more significant, and so more large finite clusters will freeze.  Super-critical PCF on $\ZZ^d$ must also pass through a critical time, and therefore it seems reasonable to expect that PCF on $\ZZ^d$ will behave in a similar way.

\begin{figure}[ht]
\centering
\includegraphics[width=0.9\textwidth]{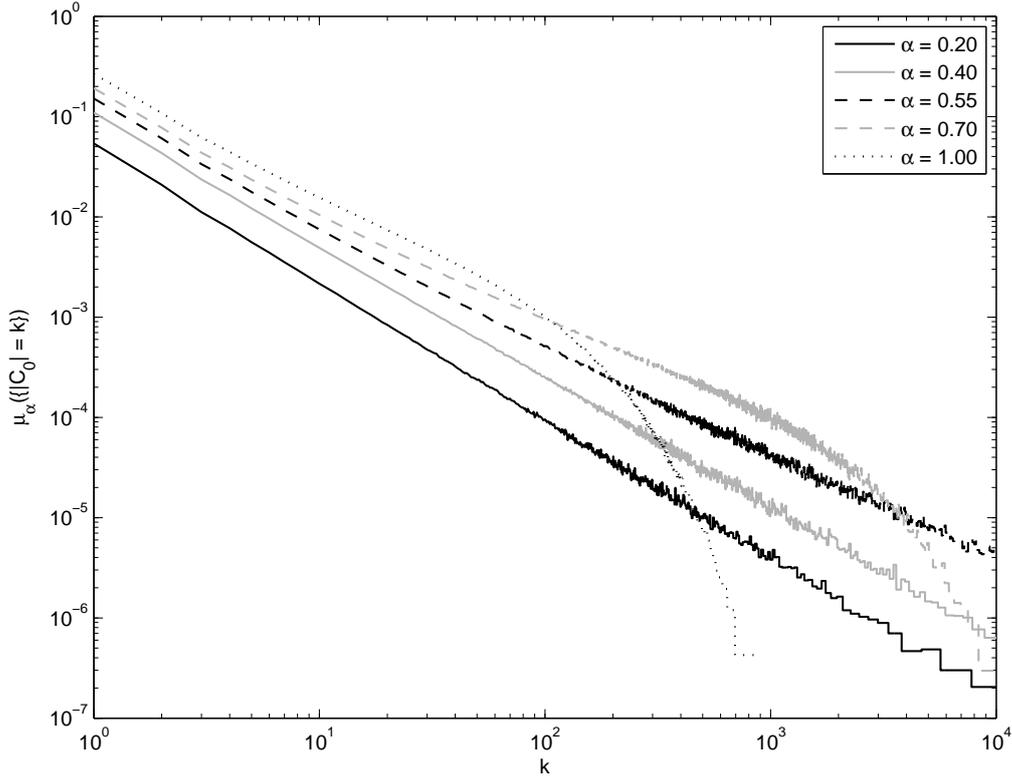}
\caption{Histogram plot with log--log axis of the cluster size distribution when running PCF on a $1024$ by $1024$ square grid at rates $\alpha=0.15$, $0.30$, $0.55$, $1.00$ and $2.00$.  Note that the data has been compiled from $1000$ simulations for each value of $\alpha$, and the widths of the bars are set so that each represents at least $100$ clusters.}\label{fig:sizes}
\end{figure}

In the case of $\ZZ^2$, the shape of the super-critical cluster of Figure \ref{fig:PCF050} (with its large voids) goes some way to showing that super-critical PCF contains large finite clusters, and so we investigate this further using Monte Carlo simulations on a finite square grid.  Figure \ref{fig:sizes} presents our results as a log--log histogram.

The curves for $\alpha=1.00$ and $\alpha=0.70$ correspond to sub-critical PCF, $\alpha=0.55$ corresponds to near critical PCF, and the values $\alpha=0.40$ and $\alpha=0.20$ correspond to super-critical PCF.  Thus we see strong numerical evidence that the cluster size distribution of PCF on $\ZZ^2$ behaves in a similar way to that of PCF on a $d$-ary tree:
\vspace{-0.5\baselineskip}\begin{itemize}
\item The plots for $\alpha=1.00$ and $\alpha=0.70$ both curve downwards corresponding to exponential decay in the size of the components.  Moreover this rate of decay is faster when the rate of freezing is faster as we would expect.
\item The log-log plot for the near-critical $\alpha=0.55$ appears to give a straight line, suggesting that the component size has a power law distribution.  Also note that the gradient is slightly less than $-1$ and thus the exponent for $\mu_{\alpha}(\{|C_0|=k\})$ is slightly larger than $1$ -- as is also the case in bond percolation.
\item The plots for $\alpha=0.20$ and $\alpha=0.40$ also follow straight lines, demonstrating that the cluster size in the super-critical case also obeys a power law distribution.  We notice further that the lines for $\alpha=0.15$ and $\alpha=0.30$ appear to have almost the same gradient as each other, but are both slightly steeper than in the near-critical regime.  This suggests that we again have a larger exponent for $\mu_{\alpha}(\{|C_0|=k\})$ in the super-critical regime than at criticality.
\end{itemize}

We now recall that the decay exponent of PCF on a $d$-ary tree at criticality is $2-\frac{1}{2d}$, since this depends on the dimension and is never constant for finite $d$ it therefore seems reasonable to assume that critical PCF on $\ZZ^d$ will never reach a mean-field value.  However, in the super-critical  regime we still have a fixed exponent (of $2$) on a tree, and therefore we would expect that in sufficiently large dimension super-critical PCF on $\ZZ^d$ will also exhibit mean field behaviour.  Our observations therefore lead us to make the following hypothesis.

\begin{hypothesis}\label{conj:sizes}
Suppose we run PCF on $\ZZ^d$ and let $C_0$ be the cluster containing the origin in the final distribution, then its size obeys
\begin{align*}
\mu_{\ZZ^d,\alpha}(\{|C_0|=k\}) \left\{ \begin{array}{ll} \prec  \ee^{-\lambda(\alpha,d)k} & \text{for } \alpha > \alpha_c  \\ \asymp  k^{-\gamma(d)} & \text{for } \alpha=\alpha_c \\ \asymp  k^{-\delta(d)} & \text{for } \alpha < \alpha_c  \end{array} \right. .
\end{align*}
Where $\gamma(d)$ never attains a mean field value but satisfies $\gamma(d)<\delta(d)\leq 2$ for all $d$, and   $\delta$ is such that $\delta(d)=2$ for $d$ sufficiently large.
\end{hypothesis}